\documentclass[a4paper,11pt]{article}

\usepackage[T1]{fontenc}
\usepackage[utf8]{inputenc}
\usepackage{lmodern}
\usepackage[english]{babel} 
\usepackage{amsfonts,amssymb,amsmath,dsfont,amsthm}
\usepackage{mathrsfs}  
\usepackage{upgreek}
\usepackage{mathtools}
\usepackage{relsize,tikz}
\usepackage{color}
\usepackage{graphicx}
\usepackage{microtype}
\usepackage{mathtools}
\usepackage[colorlinks=true, pdfstartview=FitV, linkcolor=blue, citecolor=blue, urlcolor=blue,pagebackref=false]{hyperref}
\mathtoolsset{showonlyrefs}
\usepackage{soul}

\def\<{\langle}
\def\>{\rangle}
\def\chv#1{\{\,#1\,\}}
\def\abs#1{\left\vert #1 \right\vert} 
\def\norm#1{\left\Vert #1 \right\Vert} 
\def\prt#1{\left( #1 \right)} 
\def\crt#1{\left[ #1 \right]} 
\def\set#1{\left\{\, #1 \,\right\}} 

\newcommand{\mc}[1]{{\mathcal #1}}

\newcommand{\bb}[1]{{\mathbb #1}}

\newcommand{\gep}{\varepsilon}

\newlength{\dhatheight}

\definecolor{labelkey}{gray}{.8}
\definecolor{refkey}{gray}{.8}

\definecolor{darkred}{rgb}{0.9,0.1,0.1}

\numberwithin{equation}{section}

\newtheorem{theorem}{Theorem}[section]
\newtheorem{lemma}[theorem]{Lemma}
\newtheorem{proposition}[theorem]{Proposition}
\newtheorem{corollary}[theorem]{Corollary}

\theoremstyle{definition}

\theoremstyle{remark}
\newtheorem{remark}[theorem]{Remark}

\newcommand{\vertiii}[1]{{\left\vert\kern-0.25ex\left\vert\kern-0.25ex\left\vert #1
    \right\vert\kern-0.25ex\right\vert\kern-0.25ex\right\vert}}

\def \C {{\mathbb C}}
\def \E {{\mathbb E}}

\def \P {{\mathbb P}}

\def \R {{\mathbb R}}

\def \X {\zeta}
\def \N {{\mathbb N}}
\def \V {{\mathbb V}}
\def \i {x} 
\def \j {y} 

\def \m {\mathsf{a}} 
\def \C {\mathsf{C}} %

\begin{document}

\title{Reaction-diffusion models for a class of infinite-dimensional non-linear stochastic differential equations}
\date{\today}

\author{\renewcommand{\thefootnote}{\arabic{footnote}}
Conrado da Costa\footnotemark[1],~ Bernardo Freitas Paulo da Costa\footnotemark[2],~ Daniel Valesin\footnotemark[3]}
\footnotetext[1]{
  Department of Mathematical Sciences,
  Mathematical Sciences \& Computer Science Building
  Durham University
  Upper Mountjoy Campus
  Stockton Road
  Durham University
  DH1 3LE.}
\footnotetext[2]{Instituto de Matemática - UFRJ -  Centro de Tecnologia - Bloco C
  Av. Athos da Silveira Ramos 149
  Cidade Universitária, Rio de Janeiro - RJ - Brazil.}
\footnotetext[3]{Bernoulli Institute,
  University of Groningen, Nijenborgh 9 9747 AG Groningen.}

\maketitle

\begin{abstract}
We establish the existence of solutions to a class of non-linear stochastic differential equation of reaction-diffusion type
in an infinite-dimensional space, with diffusion corresponding to a given transition kernel.
The solution obtained is the scaling limit of a sequence of interacting particle systems,
and satisfies the martingale problem corresponding to the target differential equation.


\end{abstract}

\section{Introduction}
\subsection{Background}
In this paper, we extend the main results from~\cite{dCosFreJar18} to reaction-diffusion systems evolving on infinite sets.
As in~\cite{dCosFreJar18}, the class of stochastic differential equations we consider here is:
\begin{equation}\label{eq:sde}
\begin{cases}d\zeta_t(x) = \left(\Delta_p \zeta_t(x) - b \cdot (\zeta_t(x))^\kappa \right)dt + \sqrt{a \cdot (\zeta_t(x))^\ell}\;dB^x_t,\quad x \in \V,\\[.2cm]
\zeta_0 = \bar{\zeta} \in [0,\infty)^\V,
\end{cases}
\end{equation}
where~$a,b,\kappa,\ell$ are positive real numbers
with~$\kappa,\ell \ge 1$,~$\V$ is a discrete set,~$\Delta_p$
is the Laplacian induced by a probability kernel~$p$ on~$\V$,
that is,
\begin{equation}\label{e:deltap}
  \Delta_p\zeta(x):= \sum_{y \in \V}\big(p(y,x)\cdot \zeta(y) - p(x,y)\cdot \zeta(x)\big),
 \end{equation}
and~$\{B^x_\cdot\}_{x \in \V}$ is a family of independent standard Brownian motions on~$\mathbb{R}$.

This system of equations can be used to model a reaction-diffusion
system associated for instance with chemical reactions or population dynamics.
In the setting of chemical reactions, space is divided into cells,
corresponding to points of~$\V$, and each cell contains a certain density of particles.
Within each cell, particles are subject to a reaction that can lead to a change in their density.
As an image, consider the evolution of the density of ozone subject to the
reaction ozone $\leftrightharpoons$ oxygen in a confined region.
This modelling framework is inspired by auto-catalytic models as presented
in Nicolis and Prigogine~\cite[Chap. 7]{NicPri77},
and resembles the modelling adopted by Blount~\cite{Blo96},
the main difference being that we keep the size of reaction cells constant.
Note that the system in~\eqref{eq:sde} has~$\zeta \equiv 0$ as a stable point,
and the interaction term~$-b(\zeta_t(x))^\kappa$) can then be interpreted
as a restoring force, driving the system back to equilibrium.
Hence, a solution to~\eqref{eq:sde}
represents how these processes converge to 
equilibrium in a path-wise sense.

The focus of~\cite{dCosFreJar18} was on the finite-dimensional setting,
that is, when~$\V$ is finite.
Therein a sequence of interacting particle systems~$\{\eta^n_\cdot\}_{n \ge 1}$ on~$\{0,1,\ldots\}^\V$
is shown to converge after being properly rescaled to a solution of~\eqref{eq:sde}.
This solution was moreover proved to be unique.
For each~$n$, the dynamics of~$\eta^n_\cdot$ can be encoded by the following formal generator expression,
for~$\eta \in \{0,1\ldots\}^\V$ and a local function~$f:\{0,1\ldots\}^\V \to \mathbb{R}$:
\begin{equation}\label{eq:formal_generator}
  \begin{split}
    L^nf(\eta) &= \sum_{x,y \in \V} \eta(x)\cdot p(x,y)\cdot  (f(\eta + \delta_y - \delta_x) - f(\eta)) \\
    &\quad + \sum_{x \in \V} \left[ F^{n,+}(\eta(x))\cdot (f(\eta + \delta_x) - f(\eta)) + F^{n,-}(\eta(x))\cdot (f(\eta - \delta_x) - f(\eta)) \right].
  \end{split}
\end{equation}
In words, a pile of~$\eta_t^n(x)$ particles occupies site~$x$ at time~$t$;
each particle moves with rate one according to the kernel~$p$,
and in addition, particles are born and die at~$x$ with rates~$F^{n,+}(\eta^n_t(x))$ and~$F^{n,-}(\eta^n_t(x))$, respectively.
The motion of distinct particles and births and deaths at distinct sites are independent.

The functions~$F^{n,+}$ and~$F^{n,-}$ are defined, for every~$u \ge 0$, as
\begin{align}
  \label{eq:def_F+}&F^{n,-}(u):= \frac{an^2}{2}\cdot \left(\frac{u}{n}\right)^\ell + \min\left\{\frac{an^2}{2}\cdot \left(\frac{u}{n}\right)^\ell
                     ;\;\frac{bn}{2}\cdot\left(\frac{u}{n}\right)^\kappa \right\},\\ 
  \label{eq:def_F-}&F^{n,+}(u):= \frac{an^2}{2}\cdot \left(\frac{u}{n}\right)^\ell - \min\left\{\frac{an^2}{2}\cdot \left(\frac{u}{n}\right)^\ell
                     ;\;\frac{bn}{2}\cdot\left(\frac{u}{n}\right)^\kappa \right\}.
\end{align}
These rates are chosen so that, for every~$z \ge 0$, we have
\begin{align*}
  &\lim_{n \to \infty} \frac{1}{n}(F^{n,+}(nz) - F^{n,-}(nz)) = -bz^\kappa
    \quad\text{and} \quad
    \lim_{n \to \infty} \frac{1}{n^2}(F^{n,+}(nz) + F^{n,-}(nz)) = az^\ell.
\end{align*}
The idea is to make it so that the interacting particle system resembles, with increasing precision as~$n \to \infty$, a solution to~\eqref{eq:sde}.
In addition, these functions satisfy other important properties. First,~$F^{n,-}(0) = F^{n,+}(0)= 0$,
so there is no birth (and evidently no death) of particles at empty sites.
Second,~$F^{n,-}(u) \ge F^{n,+}(u) \ge 0$ for all~$u$;
this guarantees that the number of particles in the system is stochastically decreasing, so that the dynamics has no finite-time explosion.

This leads to the result that, given a sequence~$\{\eta^n_0\}_{n \ge 1}$ with~$\frac{1}{n}\eta^n_0 \to \bar{\zeta} \in [0,\infty)^\V$,
and letting~$\eta^n_\cdot$ denote the process started from~$\eta^n_0$ and with dynamics governed by~$L^n$,
we have that the sequence of processes~$\{\tfrac{1}{n}\eta^n_\cdot\}_{n \ge 1}$
converges to a weak solution of~\eqref{eq:sde} \cite[Theorem~1]{dCosFreJar18}.
We will review the meaning of a weak solution of the SDE~\eqref{eq:sde}
in Section~\ref{s:prop_sol} below.
A limit obtained through this sort of scaling procedure,
where there is no scaling of space, but the ``mass'' of individual particles is taken to zero,
is often referred to as a \textit{fluid limit}.

\subsection{Results}

Here we are interested in obtaining the fluid limit described above
in the case where~$\V$ is a countably infinite set.
Through this extension, one can hope to achieve a better understanding of stability properties
of the solution with respect to the underlying space.
Apart from this, the extension has theoretical interest, as it brings forward some important challenges.

Our approach requires an assumption on the transition kernel~$p(\cdot,\cdot)$,
as well as a restriction on the set of allowed initial conditions for the SDE.
As in \cite{LigSpi81}, we assume that there exists a function~$\alpha:\V \to (0,\infty)$ such that
\begin{equation}
  \label{eq:main_assumption}
  \sup_{x \in \V}\alpha(x) < \infty
  \qquad\text{and}\qquad
  \C := \sup_{x\in\V}\sum_{y \in V}p(x,y)\frac{\alpha(y)}{\alpha(x)} < \infty.
\end{equation}
For instance, if~$\V = \mathbb{Z}^d$
and~$p(x,y) = \frac{1}{2d}\cdot \mathds{1}{\{x \sim y\}}$
(nearest-neighbours diffusion),
then these conditions are satisfied by~$\alpha(x) := \exp\{-|x|\}$,
where~$|\cdot\|$ denotes any norm in~$\R^d$.
Next, we define
\begin{equation}\label{e:0}
  \|\zeta\|:= \sum_{x \in \V}\alpha(x)\cdot |\zeta(x)| \in [0,\infty], \qquad \zeta \in \R^{\V}
\end{equation}
and the set of configurations
\begin{equation}\label{e:1}
  \mathcal{E} := \left\{\zeta \in [0,\infty)^{\V}:\;\sum_{x \in \V}\alpha(x)\cdot \zeta(x)<\infty\right\}.
\end{equation}
We will only consider initial conditions of~\eqref{eq:sde} belonging to~$\mathcal{E}$.
Assumptions and restrictions  of this type are common in the treatment of systems
involving diffusions on infinite environments; see for instance~\cite{And82, LigSpi81}.
The key point here is to avoid \textit{explosion from diffusion},
that is, situations where infinite amounts of mass can enter a finite set instantaneously
due to excessive growth of the initial configuration.

The next step in establishing a fluid limit result is to construct processes~$\eta^n_\cdot$ on~$\N_0^\V$
whose limit should be a solution to the stochastic differential equation.
As the dynamics linked to the generator~\eqref{eq:formal_generator}
has unbounded jump rates, and the space~$\N_0^\V$ is not compact (or locally compact),
such a  construction does not fall into the most standard framework of the theory,
by means of the Hille-Yosida theorem, as presented in Chapter~I of~\cite{LigIPS}.
While there are still ways to construct the process under our assumptions
(see for instance Chapter~IX of~\cite{LigIPS}, or the aforementioned references~\cite{And82, LigSpi81}),
here we avoid this construction issue by only constructing particle systems with finite mass. That is, we define the (countable) set
\begin{equation}\label{e:2}
E:= \left\{\eta \in \N_0^\V:\; \sum_{x \in \V}\eta(x) < \infty\right\}
\end{equation}
and only consider particle systems~$\eta^n_\cdot$ with initial configuration in~$E$.
This way, since the dynamics of~\eqref{eq:formal_generator} causes the number of particles to decrease stochastically,
it ends up producing a non-explosive continuous-time Markov chain on the countable state space~$E$. 

We are now ready to state our main result.

\begin{theorem}
\label{thm:main}\ 
\begin{itemize}
\item[(a)]
Let~$\zeta_0 \in \mathcal{E}$ and let~$\{\eta^n_0\}$ be a sequence in~$E$
with~$\|\tfrac{1}{n}\eta^n_0 - \zeta_0\| \xrightarrow{n \to \infty} 0$.
For each~$n$, let~$(\eta^n_t)_{t \ge 0}$ denote the Markov chain on~$E$
with transitions encoded by~\eqref{eq:formal_generator} started from~$\eta^n_0$.
Then, as~$n \to \infty$, the processes~$\eta^n_\cdot$ converge in distribution
(with respect to the Skorokhod topology) to an~$\mathcal{E}$-valued process~$\zeta_\cdot$
with continuous trajectories which is a weak solution to~\eqref{eq:sde}
with initial condition~$\bar{\zeta} = \zeta_0$. The law of this process
does not depend on the choice of sequence~$\{\eta^n_0\}$
with~$\|\tfrac{1}{n}\eta^n_0 -\zeta_0\|\xrightarrow{n \to \infty} 0$.
\item[(b)] In case~$\sum_{x \in \V}\zeta_0(x) < \infty$, the process
obtained through this limit is the unique weak solution
to~\eqref{eq:sde} with initial condition~$\bar{\zeta} = \zeta_0$
in the sense that it has the same distribution as any other solution of the same equation.
\item[(c)] The mapping~$\mathcal{E}\ni\zeta_0 \mapsto (\zeta_t)_{t \ge 0}$
of initial conditions to corresponding solutions obtained through the limit of part~(a) is continuous
when~$\mathcal{E}$ is endowed with the norm~$\|\cdot\|$ and the set of processes on~$C([0,\infty),\mathcal{E})$
is endowed with the topology of weak convergence of probability measures.
\end{itemize}
\end{theorem}

\subsection{Outline of methods and organization of the paper}
Let us give an outline of our methods. Using generator estimates,
we prove that a collection of processes~$\{\eta^n_\cdot\}$ as in the statement
of Theorem~\ref{thm:main}(a),
with~$\|\tfrac{1}{n}\eta^n_0 - \zeta_0\| \xrightarrow{n \to \infty} 0$
for some~$\zeta_0 \in \mathcal{E}$, is tight.
This allows us to extract convergent subsequences,
say~$\{\eta^{n_k}_\cdot\}$ converging to a process~$\zeta_\cdot$.
We then prove that the Dynkin martingales associated
to~$\eta^{n_k}_\cdot$, as defined in Lemma~\ref{lem:first_pre_mcm},
converge to processes of the form
\begin{equation}\label{eq:limit_dynkin}
f(\zeta_t)-f(\zeta_0)-\int_0^t \mathcal{L}^*f(\zeta_s)ds,\qquad t \ge 0,
\end{equation}
where~$\mathcal{L}^*$, the generator associated to~\eqref{eq:sde}, is given by
\begin{equation}\label{eq:want_gen} 
  (\mathcal{L}^*f)(\zeta) := \sum_{x \in \V} (\Delta_p\zeta(x) - b\cdot (\zeta(x))^\kappa)\cdot \partial_xf(\zeta)
  + \frac12\sum_{x \in \V} a \cdot (\zeta(x))^\ell\cdot \partial^2_x f(\zeta),\quad \zeta \in \mathcal{E},
\end{equation}  
for a suitable collection of functions~$f$.
This convergence allows us to obtain that~\eqref{eq:limit_dynkin} is a local martingale.
Using classical results from the theory of stochastic differential equations,
we then conclude that the subsequential limit~$\zeta_\cdot$ is a solution of~\eqref{eq:sde}.
An adaptation of the argument in~\cite{WatYam71} gives us Theorem~\ref{thm:main}(b),
that is, that~\eqref{eq:sde} has at most one solution in case~$\zeta_0$ has finite mass.
Combining these ideas, we get that if~$\zeta_0$ has finite mass,
then~$\{\eta^n_\cdot\}$ has a single accumulation point, so the whole sequence converges.
From this, we finish the proof of Theorem~\ref{thm:main}(a), that is,
we prove convergence for any $\zeta_0 \in \mathcal{E}$ with infinite mass by approximation%
: any $\zeta_0 \in \mathcal{E}$ is arbitrarily close to configurations with finite mass.

A key tool that we rely on for this approximation and for several other arguments is a coupling inequality,
Lemma~\ref{lem:sup0} below, allowing us to compare pairs
of processes with same generator but different initial configurations.

The rest of the paper is organized as follows.
In Section~\ref{s:tech}, we review several technical concepts and results,
including notions of convergence of probability measures,
local martingales defined from Markov chains,
and classical results about stochastic differential equations. In Section~\ref{s:ips}
we study particle systems with finite mass on~$\N_0^\V$,
and obtain the key coupling inequality in Lemma~\ref{lem:sup0}.
In Section~\ref{ss:sequence} we state our tightness result and
use it to follow the rest of the outline given above, proving our main results.
In Section~\ref{s:tight} we prove the tightness result. Section~\ref{s:appendix}
is an appendix where we include some proofs to ease the flow of the exposition in the paper.

\section{Technical preliminaries}\label{s:tech}
In this Section, we collect remarks, definitions, and properties
that will be useful in the study of convergence
of a family of stochastic processes as mentioned in the previous Section.

\subsection{Configuration spaces}
We let~$\V$ be a countable set
and~$p: \V \times \V \to [0,1]$ be a probability transition function
(that is, $\sum_y p(x,y) = 1$ for all~$x$),
and assume that there exists a function~$\alpha: \V \to [0,\infty)$
for which~\eqref{eq:main_assumption} holds.
We define $\|\cdot\|$, $\mathcal{E}$ and~$E$
as in~\eqref{e:0}--\eqref{e:2}.
Note that~$E$ is countable,
that~$\|\cdot\|$ is a norm on the linear subspace of~$\R^\V$ where it is finite,
and that the metric induced by~$\|\cdot\|$
turns~$\mathcal{E}$ into a complete and separable metric space.
For the sake of clarity, we mostly denote the (integer-valued) elements of~$E$ by the letter~$\eta$ rather than~$\zeta$,
and processes taking values on~$E$ by~$\eta_\cdot$ rather than~$\zeta_\cdot$.

It will be useful to observe that the assumptions~\eqref{eq:main_assumption} yields:
\begin{equation}\label{eq:quick_bound}
p(x,y) \le \frac{1}{\alpha(y)}\sum_{z} p(x,z)\cdot \alpha(z)\le \frac{\C \alpha(x)}{\alpha(y)}.
\end{equation}
This implies that, for~$\zeta\in \mathcal{E}$,
\[
  \sum_{x \in \V} \zeta(x) \cdot p(x,y)  \le
  \C\cdot \sum_{x \in \V} \frac{\zeta(x)\cdot \alpha(x)}{\alpha(y)}= \C \frac{\|\zeta\|}{\alpha(y)},
  \qquad y \in \V.
\]
In particular, $\Delta_p \zeta(x)$  in~\eqref{e:deltap}
is well defined for all~$\zeta \in \mathcal{E}$ and all~$x \in \V$.
\subsection{Convergence of probability measures on trajectory spaces}
To study convergence of probability measures on trajectory spaces,
we first define a metric on the space of trajectories,
then we consider a family of $\sigma$-algebras associated to this metric
and finally we define a distance between probability measures on such $\sigma$-algebras.

\paragraph{Metric.}
Let $\mathcal{X} = (\mathcal{X},\mathrm{d}_\mathcal{X})$
be a complete, separable metric space.
In most cases,
this will be either~$(\R,|\cdot|)$
or~$\mathcal{E}$ or~$E$ with the metric induced by~$\|\cdot\|$.
We denote by~$D_{\mathcal{X}} = D([0,\infty), \mathcal{X})$
the space of càdlàg functions $\gamma : [0,\infty) \to \mathcal{X}$,
and by~$C_{\mathcal{X}}$ the set of functions in~$D_{\mathcal{X}}$
which are continuous.
The Skorokhod metric on~$D_{\mathcal{X}}$ is defined by
\[
\mathrm{d}_{\mathrm{S}}(\gamma, \gamma'):= \int_0^\infty e^{-t}\cdot \mathrm{d}^{(t)}_{\mathrm{S}}(\gamma,\gamma')dt,
\]
where
\[
  \mathrm{d}^{(t)}_{\mathrm{S}}(\gamma,\gamma'):=
  1 \wedge \inf_\varphi \left( \sup_{s \in [0,t]}\mathrm{d}_{\mathcal{X}}( \gamma_{\varphi(s)}, \gamma'_{s})
    \vee \sup_{r,s \in [0,t]} \log\frac{|\varphi(r)-\varphi(s)|}{|r-s|} \right),
\]
where the infimum is taken over all increasing bijections~$\varphi:[0,t] \to [0,t]$.
This turns $(D_{\mathcal{X}},\mathrm{d}_{\mathrm{S}})$
into a complete and separable metric space,
and we denote by $\mc{D}_{\mathcal{X}}$ its Borel $\sigma$-algebra.
We refer the reader to~\cite[Chapter 3]{Bil99} and~\cite[Chapter 3]{EthKur86} for expositions on this metric.
Here let us only make one further observation, see~\cite[Section 12, p. 124]{Bil99}:
\begin{multline}
  \label{eq:ds_cont}
  \text{if } \{\gamma^n\}_{n \ge 1} \subset D_{\mathcal{X}},\; \gamma \in C_{\mathcal{X}}
  \text{ and } \mathrm{d}_{\mathrm{S}}(\gamma^n,\gamma) \xrightarrow{n \to \infty}0, \\[.2cm]
  \text{then } \sup_{0 \le s \le t} \mathrm{d}_{\mathcal{X}}(\gamma^n_s,\gamma_s)
    \xrightarrow{n \to \infty} 0\; \text{for all } t\ge 0,
\end{multline}
that is, convergence in the Skorokhod topology to a continuous function implies uniform convergence on compact intervals.

\paragraph{Sigma-algebras.}
Given a stochastic process $X_\cdot$ on $\mc{X}$ with càdlàg trajectories,
we denote by $\mc{F}_t = \mc{F}^X_t$ the $\sigma$-algebra
generated by $(X_s)_{0\leq s \leq t}$
and by~$\mc{N} := \{A \in \mc{D}_{\mathcal{X}} \colon \bb{P}(X_\cdot \in A) = 0\}$.
We refer to $(\mc{F}_t)_{t \geq 0}$ as
the \emph{natural filtration} of $(X_t)_{t \geq 0}$.

\paragraph{Convergence.}
Finally, we recall the definition of the Lévy-Prohorov distance
in the case of two probability measures $\mu$ and~$\nu$
defined on $D_\mc{X} = ({D}_{\mathcal{X}},\mc{D}_{\mathcal{X}})$:
\[
  \mathrm{d}_{\mathrm{LP}}(\mu, \nu)
  := \inf\left\{\varepsilon > 0 : \; \mu(A) \le \nu(A^\epsilon) + \epsilon \text{ and } \nu(A) \le \mu(A^\epsilon) + \epsilon
                \text{ for all } A \in \mathcal{D}_{\mathcal{X}}\right\},
\]
where~$A^\epsilon:= \{ y\in D_{\mc{X}} :\mathrm{d}_{\mathrm{S}}(x,y) < \epsilon \text{ for some }x \in A\}$.
Convergence in this metric is equivalent to
weak convergence of probability measures, see~\cite[Theorem 3.3.1, p. 108]{EthKur86},
that is,~$\mathrm{d}_{\mathrm{LP}}(\mu_n,\mu)\xrightarrow{n \to \infty} 0$
is equivalent to having~$\int f\,\mathrm{d}\mu_n \xrightarrow{n \to \infty} \int f\, \mathrm{d}\mu$
for all continuous and bounded functions~$f:D_{\mathcal{X}} \to \mathbb{R}$.
Denote by $\mc{C}(\mu,\nu)$ the set of all measures $\hat{\lambda}$ on $D_{\mathcal{X}} \times D_{\mathcal{X}}$ that couple $\mu,\nu$. By~\cite[Theorem 3.1.2, p. 98]{EthKur86}, we remark that
  \begin{equation}\label{LPcoup}
    \mathrm{d}_{\mathrm{LP}}(\mu, \nu) = \inf_{\lambda \in \mc{C}(\mu,\nu)}
    \quad\inf\big\{\varepsilon>0 \colon \lambda\{(\gamma,\gamma') \in D_{\mathcal{X}} \times D_{\mathcal{X}}:
    \mathrm{d}_{\mathrm{S}}(\gamma, \gamma')\geq \varepsilon\} \leq \varepsilon\big\}.  
  \end{equation}

\subsection{Continuous-time Markov chains and martingales}
Given a stochastic process $X$ on $D_{\mc{X}}$ there is a sequence of stopping times
$\tau^X_0 := 0 $ and $\tau^X_n := \inf\{t > \tau^X_{n-1} \colon X_t \neq X_{t-}\}$
that exhaust the jumps of $X$, see Proposition 2.26 in~\cite[Chapter 1, p. 10]{KarShr98}.
Let $\tau^X_\infty := \lim_n \tau^X_n$.
Following~\cite[Remark 2.27]{Lig10}, we say that a process on $D_{\mc{X}}$ is \emph{non-explosive} if
\begin{equation}\label{nonex}
  \bb{P}(\tau^X_\infty < t ) = 0\quad \text{ for all } t \geq 0.
\end{equation}
    
For the following two results, let~$S$ be a countable set and~$(X_t)_{t \ge 0}$ be a non-explosive continuous-time Markov chain on~$S$.
For distinct~$x,y \in S$, let~$q(x,y) \ge 0$ be the jump rate from~$x$ to~$y$, and let~$q(x,x) = -\sum_{y \neq x}q(x,y)$.
For a function~$f: S \to \mathbb{R}$ satisfying
\begin{equation}\label{eq:conv_gener}
\sum_{y \in S}q(x,y)\cdot |f(y)-f(x)| < \infty\qquad \text{for all }x \in S,
\end{equation}
we define $Lf(x):= \sum_{y\in S} q(x,y)\cdot (f(y)-f(x))$.
If~$f$ also satisfies
\begin{equation}\label{eq:conv_quad}
\sum_{y \in S}q(x,y)\cdot (f(y)-f(x))^2 < \infty\qquad \text{for all }x \in S,
\end{equation}
we define~$Qf(x):= \sum_{y\in S} q(x,y)\cdot (f(y)-f(x))^2$.

\begin{lemma}\label{lem:first_pre_mcm}
Let~$f: S \to \mathbb{R}$ be a function satisfying~\eqref{eq:conv_gener}.
Then, the process
\[
M^f_t:= f(X_t) - f(X_0) - \int_0^t Lf(X_s)ds,\qquad t\ge 0
\]
is a local martingale with respect to the natural filtration of~$(X_t)_{t \ge 0}$.
If~$f$ also satisfies~\eqref{eq:conv_quad}, then the process
\[
N^f_t:= M_t^2 - \int_0^t Qf(X_s)ds,\qquad t \ge 0
\]
is also a local martingale with respect to the natural filtration of~$(X_t)_{t \ge 0}$.
\end{lemma}
\begin{proof}
Fix an arbitrary initial state~$x_0 \in S$,
and let~$(\Lambda^j)_{j \ge 1}$ be an increasing sequence of finite subsets of~$S$
with~$x_0 \in \Lambda^1$ and~$\cup_j\Lambda^j = S$.
Define~$\tau^j := \inf\{t \ge 0:\;X_t \notin \Lambda^j\}$.
We then have that~$\tau^j \le \tau^{j+1}$ for each~$j$ and,
because $(X_t)_{t\geq 0}$ is  non-explosive,
$\tau^j \xrightarrow{j \to \infty} \infty$ almost surely.
Then, under the assumptions~\eqref{eq:conv_gener} and~\eqref{eq:conv_quad},
classical arguments establish that~$M^f_\cdot$ and~$N^f_\cdot$
are local martingales  with respect to the natural filtration of $(X_t)_{t \geq 0}$,
see for instance~\cite[Appendix 1, Lemma 5.1]{KipLan99}.
\end{proof}
\begin{remark}
  The martingales $M^f_\cdot$ for functions $f$ that satisfy~\eqref{eq:conv_gener}
  are commonly referred to as \emph{Dynkin Martingales}.
  Unless the generator $L$ is not clear from the context,
  we will omit the dependence (as we have done here) to alleviate notation.
\end{remark}

To obtain stochastic bounds,
we prove a supermartingale inequality associated to the martingales $M^f_\cdot$.
\begin{lemma}\label{lem:super_from_mc}
Let~$f: S \to \mathbb{R}$ be a non-negative function satisfying~\eqref{eq:conv_gener}.
Assume that there exists~$C > 0$ such that
\[
\sum_{y \in S}q(x,y)\cdot (f(y)-f(x)) \le Cf(x) \qquad \text{ for all }x \in S.
\]
Then, the process~$(e^{-Ct}\cdot f(X_t))_{t \ge 0}$ is a supermartingale with respect to the natural filtration of~$(X_t)_{t \ge 0}$.
\end{lemma}
\begin{proof}
Fix an arbitrary initial state~$x_0 \in S$.
By the Markov property, it is sufficient to prove that, for any~$t \ge 0$,
\begin{equation}
\mathbb{E}[e^{-Ct}\cdot f(X_t)] \le f(x_0).
\end{equation}
Let~$(\Lambda^j)_{j \ge 1}$ and~$\tau^j$ be as in the proof of Lemma~\ref{lem:first_pre_mcm}. Define the processes
\[
  X^j_t :=
  \begin{cases}
    X_t&\text{if }t \le \tau^j;\\
    \star&\text{otherwise,}
  \end{cases}
  \qquad t\ge 0
\]
where~$\star$ denotes a cemetery state.
Then,~$X^j_\cdot$ is a Markov chain on the finite state space~$S^j := \Lambda^j \cup \{\star\}$,
with jump rates, for distinct~$x,y \in S_j$, given by
\[
  q^j(x,y) =
  \begin{cases}
    0&\text{if }x = \star;\\
    q(x,y)& \text{if } x \neq \star,\;y \neq \star;\\
    \sum_{y \in (\Lambda^j)^c}q(x,y)&\text{if } x \neq \star,\;y = \star.
  \end{cases}
\]
Let~$f^j: S^j \to \mathbb{R}$ be defined by~$f^j(x) =  f(x)$ for~$x \neq \star$ and~$f^j(\star) = 0$.
We have, for any~$x \in S^j\backslash \{\star\}$,
\[
  \sum_{y \in S^j} q^j(x,y)\cdot (f^j(y)-f^j(x)) \le \sum_{y \in S} q(x,y)\cdot (f(y)-f(x)) \le Cf(x),
\]
since~$0 = f^j(\star) \le f(y)$ for~$y \in S \backslash \Lambda^j$.
It then follows from the elementary theory of Markov chains that~$(e^{-Ct}\cdot X^j_t)_{t \ge  0}$
is a supermartingale. Since~$f(X^j_t)\xrightarrow{j \to \infty} f(X_t)$ almost surely,  it follows from Fatou's Lemma that
\[
  \mathbb{E}[e^{-Ct}\cdot f(X_t)] \le \liminf_{j \to \infty} \mathbb{E}[e^{-Ct}\cdot f(X_t^j)] \le f(x_0).
\]
\end{proof}

\subsection{Some properties of solutions of the SDE (\ref{eq:sde})}\label{s:prop_sol}

Let~$(\zeta_t)_{t \ge 0}$ be a stochastic process (defined on some space~$(\Omega,\mathcal{F},\P)$)
with values in~$\mathcal{E}$ and continuous trajectories. We say that~$\zeta_\cdot$
is a weak solution to the SDE~\eqref{eq:sde}
if there exists a space~$(\tilde\Omega,\tilde{\mathcal{F}},\tilde\P)$ in which we have defined
\begin{itemize}
\item a process~$X_\cdot$ with values in~$\mathcal{E}$, continuous trajectories, and same distribution as~$\zeta_\cdot$, and
\item a family~$(B^x_\cdot)_{x \in \V}$ of independent,  standard one-dimensional  Brownian motions, 
\end{itemize}
and moreover,~$\tilde \P$ almost surely we have
\[
  X_t(x) = X_0(x) + \int_0^t (\Delta_pX_s(x) - b(X_s(x))^\kappa)ds
                  + \int_0^t \left(a(X_s(x))^\ell\right)^{\frac12}dB^x_s,
  \quad \text{for $t \ge 0$ and $x \in \V$}.
\]

We will need the following result, which gives uniqueness for solutions of~\eqref{eq:sde} with finite mass.
The essence of its proof is taken from~\cite{WatYam71}.
We present in the Appendix (Section~\ref{s:unique}) the proof with slight modifications to adjust to our setting.
\begin{proposition} \label{prop:unique}
  Let~$\bar{\zeta} \in \mathcal{E}$ satisfy~$\|\zeta\|_1 = \sum_{x \in \V}\zeta(x) < \infty$. Let~$\zeta^1_\cdot$ and~$\zeta^2_\cdot$
  be two weak solutions of the SDE~\eqref{eq:sde} with initial condition~$\bar{\zeta}$.
  Then,~$\zeta^1_\cdot$ and~$\zeta^2_\cdot$ have the same distribution.
\end{proposition}

Let~$f:\mathcal{E} \to \mathbb{R}$ be a function
for which there exists a finite set~$\{x_1,\ldots, x_k\} \subset \V$ so that~$f$
only depends on~$(\zeta(x_1),\ldots, \zeta(x_k))$,
and moreover~$f$ is a twice continuously differentiable function of this vector. Define
\begin{equation} \label{eq:def_L*}
  (\mathcal{L}^*f)(\zeta) := \sum_{x \in \V} (\Delta_p\zeta(x) - b\cdot (\zeta(x))^\kappa)\cdot \partial_xf(\zeta)
  + \frac12\sum_{x \in \V} a \cdot (\zeta(x))^\ell\cdot \partial^2_x f(\zeta),\quad \zeta \in \mathcal{E}.
\end{equation}
In particular, we denote by $f_x(\zeta) := \zeta(x)$ the coordinate projection on $x \in \V$,
and by $f_{xy}(\zeta) := \zeta(x)\zeta(y)$.
Note that all those functions are indeed differentiable and depend only on a finite number of sites.

\begin{proposition}  \label{prop:are_sols}
  Let~$\zeta_\cdot$ be a stochastic  process with values in~$\mathcal{E}$ and continuous trajectories.
  Assume that for every~$f \in \{f_x,f_{xy}:x,y \in \V\}$, the process
  \[
    f(\zeta_t) - f(\zeta_0) - \int_0^t \mathcal{L}^*f(\zeta_s)ds,\qquad  t \ge 0.
  \]
  is a local martingale. Then,~$\zeta_\cdot$ is a weak solution of~\eqref{eq:sde}.
\end{proposition}
This proposition is the same as Proposition~4.6, page 315 in~\cite{KarShr98},
except that here we deal with infinite-dimensional processes,
whereas the setup of the proposition in~\cite{KarShr98} is finite-dimensional.
However, largely due to the fact that the cross-variation of our SDE is trivial
(that is, the expression for~$d\zeta_t(x)$ in~\eqref{eq:sde} does not involve~$dB^y_t$ for~$y \neq x$),
the proof in~\cite{KarShr98} carries through to our setting.
In order to highlight the differences and the main steps,
we sketch the proof in the appendix (Section~\ref{s:are_sols}).

\section{Construction of diffusive birth-and-death particle systems}
\label{s:ips}
Recall that~$E = \{\eta \in \N_0^\V:\sum_x \eta(x) < \infty\}$.
In this section, we will construct continuous-time Markov chains on~$E$
that will later be used in the fluid limit for the proof of our main result.
Rather than having an index~$n$ and functions~$F^{n,+}$ and~$F^{n,-}$ as in~\eqref{eq:def_F+} and~\eqref{eq:def_F-},
for now we will have no index, and functions~$F^+$ and~$F^-$ satisfying certain properties (see~\eqref{decreasingF} below).

We define the set of \textit{marks}
\begin{equation}\label{marks}
\mathcal{M} = \{(\i,+),(\i,-),(\i,\j)\colon \i,  \j\in \V\}.
\end{equation}
Marks will serve as instructions for the dynamics. Marks of the form~$(x,+)$ and~$(x,-)$
represent the birth and the death of a particle at site~$x$, respectively, and a mark of the form~$(x,y)$
represents that a particle from site~$x$ jumps to site~$y$. Marks are thus associated with the \textit{transition operators}
\begin{equation}\label{transition}
\Gamma^{\i,+}(\eta) :=\eta + \delta_\i,
\qquad \Gamma^{\i,-}(\eta) := \eta - \delta_\i,
\qquad \text{and}
\qquad \Gamma^{(\i,\j)}_\i(\eta) :=\eta - \delta_\i + \delta_{\j},
\end{equation}
where for $\i \in \V$, $\delta_x \in E$
is the configuration with only one particle at $\i$.
For a configuration $\eta \in E$ and $\i,\j \in \V$,
we define \textit{transition rates} by
\begin{equation}\label{rates}
  R^{\i,+}(\eta) := F^+(\eta(\i)),
  \qquad R^{x,-}(\eta) := F^-(\eta(\i)),
  \qquad R^{(\i,\j)}(\eta) := p(\i,\j) \cdot\eta(\i),
\end{equation}
where we assume that the \textit{reaction functions} $F^+$ and $F^-$ satisfy
\begin{equation}
  F^+(0) = F^-(0) = 0,\qquad  0 \leq F^+ \leq F^-,\qquad
  \N_0 \ni z \mapsto F^+(z) - F^-(z) \text{ is decreasing.} \label{decreasingF}
\end{equation}

We now define a continuous-time Markov chain on~$E$ with the prescription that
\[
\text{for each }\m \in \mathcal{M},\; \text{$\eta$ jumps to $\Gamma^\m(\eta)$ with rate $R^\m(\eta)$.}
\]
Noting that~$\sum_{\m \in \mathcal{M}}R^\m(\eta) < \infty$ for any~$\eta \in E$,
this indeed describes jump rates of a continuous-time Markov chain.
The assumption $F^+ \leq F^-$ combined with a simple stochastic comparison argument
(see~\cite[Lemma 2]{dCosFreJar18}) implies that the chain is non-explosive.
It will be important to leave the initial condition explicit,
so we will denote the chain started from~$\eta \in E$ by~$(\Phi_t(\eta))_{t \ge 0}$.

For any function~$f:E \to \mathbb{R}$ satisfying
\begin{equation}\label{eq:def_of_conv_of_L0}
  \sum_{\m \in \mathcal{M}} R^\m(\eta)\cdot \max\left(|f(\Gamma^\m(\eta))-f(\eta)|,\;|f(\Gamma^\m(\eta))-f(\eta)|^2\right) < \infty
  \quad
  \text{for all }\eta \in E, 
\end{equation}
we define
\begin{equation}
  Lf(\eta) := \sum_{\m \in \mathcal{M}} R^\m(\eta)\cdot (f(\Gamma^\m(\eta))-f(\eta)),
  \qquad
  Qf(\eta):= \sum_{\m \in \mathcal{M}} R^\m(\eta)\cdot (f(\Gamma^\m(\eta))-f(\eta))^2.
\end{equation}
Recall that a function defined on a subset of~$\mathbb{R}^\V$ is called \textit{local}
if there exists a finite set~$\V' \subset \V$ such that the function depends on~$\eta \in \mathbb{R}^\V$
only through~$(\eta(x):x \in \V')$. We then have
\begin{lemma}\label{lem:two_martingales}
Any local function~$f: E \to \mathbb{R}$ satisfies~\eqref{eq:def_of_conv_of_L0}, and the processes
\begin{align*}
  M^f_t & := f(\Phi_t(\eta)) - f(\eta) - \int_0^t Lf(\Phi_s(\eta))ds,
                \quad t \ge 0,\text{ and }\\
  N^f_t & := \left( M^f_t \right)^2 - \int_0^t Qf(\Phi_s(\eta))ds,
                \quad t \ge 0
\end{align*}
are local martingales with respect to the natural filtration of $(\Phi_t(\eta))_{t \geq 0}$.
\end{lemma}
\begin{proof}
The first statement is straightforward to check: since~$f$ is local and its argument is an element of~$E$,
the sum in~\eqref{eq:def_of_conv_of_L0} only has finitely many non-zero terms.
The second statement then follows from Lemma~\ref{lem:first_pre_mcm}.
\end{proof}
 
\begin{lemma}[1-norm bound in~$E$]
  For any~$\eta \in E$, the process~$\|\Phi_\cdot(\eta)\|_1$ is a supermartingale.
  In particular, for any~$T \ge 0$ and~$A > 0$,
\begin{equation}
\label{eq:l1norm}
\P\left(\sup_{0\le t \le T} \|\Phi_t(\eta)\|_1 > A \right) \le \frac{\|\eta\|_1}{A}.
\end{equation}
\end{lemma}  
\begin{proof}
Let~$f:E \to \mathbb{R}$ be given by~$f(\eta) := \|\eta\|_1$.
It is readily seen that~$f$ satisfies~\eqref{eq:def_of_conv_of_L0},
and then it follows from Lemma~\ref{lem:super_from_mc} that~$\|\Phi_\cdot(\eta)\|_1$ is a supermartingale.
Equation~\eqref{eq:l1norm} then follows from the optional stopping theorem.
\end{proof}

Now, given a pair~$\eta,\eta' \in E$ we will construct a coupled process
\[
  \hat{\Phi}_t(\eta,\eta') = (\hat{\Phi}_{t,1}(\eta,\eta'),\hat{\Phi}_{t,2}(\eta,\eta')),\qquad t \ge 0
\]
on~$E \times E$ so that the coordinate processes are distributed as~$\Phi_\cdot(\eta)$ and~$\Phi_\cdot(\eta')$, respectively.
The coupling is given as the Markov chain on~$E \times E$ with transition rates described by
\[
  \text{for each } \m \in \mathcal{M},\; (\eta,\eta') \text{ jumps to }
  \begin{cases} 
    (\Gamma^\m(\eta),\Gamma^\m(\eta'))&\text{with rate } \min(R^\m(\eta),R^\m(\eta'));\\
    (\Gamma^\m(\eta),\eta')&\text{with rate } \max(R^\m(\eta)-R^\m(\eta'),0);\\
    (\eta,\Gamma^\m(\eta'))&\text{with rate } \max(R^\m(\eta')-R^\m(\eta),0).
  \end{cases}
\]
We then define
\begin{equation}\label{eq:three}
  \begin{split}
\hat{L}g(\eta,\eta'):= \sum_{\m \in \mathcal{M}} &[\min(R^\m(\eta),R^\m(\eta'))\cdot (g(\Gamma^\m(\eta),\Gamma^\m(\eta')) - g(\eta,\eta') ) \\[-.4cm]
&\;\;+ \max(R^\m(\eta)-R^\m(\eta'),0)\cdot (g(\Gamma^\m(\eta),\eta') - g(\eta,\eta'))\\
&\;\;+ \max(R^\m(\eta')-R^\m(\eta),0)\cdot (g(\eta,\Gamma^\m(\eta')) - g(\eta,\eta')) ]
\end{split}
\end{equation}
for functions~$g:E\times E \to \mathbb{R}$ for which the sum on the right-hand side is absolutely convergent for all~$(\eta,\eta')$.

\begin{lemma}\label{lem:bound_couple}
For~$g(\eta,\eta'):= \|\eta - \eta'\|$, we have that~$\hat{L}g$ is well defined and satisfies
\begin{equation}
\hat{L}g(\eta,\eta') \le \C\cdot g(\eta,\eta')\quad \text{for all }\eta,\eta' \in E,
\end{equation}
where~$\C$ is the constant from~\eqref{eq:main_assumption}.
\end{lemma}
\begin{proof}
For this choice of~$g$, the first line in the right-hand side of~\eqref{eq:three} vanishes,
so we can write~$\hat{L}g(\eta,\eta') = \sum_{\m \in \mathcal{M}} (\Xi^\m_1 + \Xi^\m_2)(\eta,\eta')$,
where
\begin{align*}
   \Xi^\m_{1}(\eta,\eta') & = \left(R^\m(\eta) - \min\{R^\m(\eta),R^\m(\eta')\}\right) \cdot \left(\|\Gamma^\m(\eta) - \eta'\| - \|\eta - \eta'\| \right),\\[.2cm]
  \Xi^\m_{2}(\eta,\eta') & = \left(R^\m(\eta') - \min\{R^\m(\eta),R^\m(\eta')\}\right) \cdot \left(\|\eta - \Gamma^\m(\eta')\| - \|\eta - \eta'\| \right).
\end{align*}

We first deal with the reaction terms, that is,
the terms corresponding to marks of the form~$(x,+)$ and~$(x,-)$.
Due to the nature of the rates we have chosen,
the only cases we have to look at are those for which $\eta(x) \neq \eta'(x)$.
If~$\eta(x) > \eta'(x)$,
 the contribution from $(x,+)$ marks is:
\begin{align*}
  \big( \Xi_1^{(x,+)} + \Xi_2^{(x,+)} \big)(\eta, \eta')
  & = [F^+(\eta(x)) - F^+(\eta'(x))]_+ \cdot \alpha(x)
      + [F^+(\eta'(x)) - F^+(\eta(x))]_+ \cdot (-\alpha(x)) \\
  & = [F^+(\eta(x)) - F^+(\eta'(x))] \cdot \alpha(x).
\end{align*}
Doing similarly for marks $(x,-)$, we get:
\begin{equation}\label{smartF}
  \sum_{\sigma \in \{+,-\}} (\Xi_1^{(x,\sigma)} + \Xi_2^{(x,\sigma)})(\eta, \eta')
  = \alpha(x)\cdot \crt{(F^+-F^-)(\eta(x)) - (F^+ - F^-)(\eta'(x))} \le 0,
\end{equation}
where the last inequality follows from
our hypothesis~\eqref{decreasingF}
ensuring that $F = F^+ - F^-$ is decreasing.
Observe that we don't need to assume (although it would be natural to)
that each $F^+$ and $F^-$ are increasing.

The same argument shows \eqref{smartF} for the case~$\eta(x) < \eta'(x)$.
We thus conclude that
\begin{equation}\label{boundxt}
  \sum_{x \in \V} \sum_{\sigma \in \{+,-\}} (\Xi^{(x,\sigma)}_1 + \Xi^{(x,\sigma)}_2)(\eta,\eta') \le 0.
\end{equation}

We now turn to the diffusion terms.
Fix~$x,y \in \V$, and first assume that~$\eta(x) > \eta'(x)$.
Again, note that reaction rates are equal when $\eta(x) = \eta'(x)$,
and the contribution to $\hat{L}g$ is zero in those cases.
We then have~$\Xi^{(x,y)}_2(\eta,\eta') = 0$ and
\begin{align*}
  \Xi^{(x,y)}_1(\eta,\eta')
  & = (\eta(x)-\eta'(x)) \cdot p(x,y)
      \cdot \left(\|\eta -\delta_x+\delta_y - \eta' \| - \|\eta - \eta' \|\right) \\
  & \le (\eta(x)-\eta'(x)) \cdot p(x,y) \cdot (- \alpha(x) + \alpha(y)).
\end{align*}
Treating the case~$\eta(x) < \eta'(x)$ analogously, we obtain
\begin{align*}
  \sum_{x,y \in \V} (\Xi^{(x,y)}_1 + \Xi^{(x,y)}_2)(\eta,\eta')
  & \le \sum_{x:\eta(x) > \eta'(x)} \sum_y \Xi_1^{(x,y)}(\eta,\eta')
    \quad + \sum_{x:\eta(x) < \eta'(x)} \sum_y \Xi_2^{(x,y)}(\eta,\eta') \\[.2cm]
  & \le \sum_{x \in \V} |\eta(x) - \eta'(x)| \sum_{y \in \V} p(x,y) ( - \alpha(x) + \alpha(y)) \\
  & \stackrel{\eqref{eq:main_assumption}}{\le}
    (\C - 1)\sum_{x \in \V} |\eta(x)-\eta'(x)|\alpha(x)
    \leq \C \|\eta-\eta'\|.
  \qedhere
\end{align*}
\end{proof}
\begin{lemma}\label{lem:sup0}
For any~$\eta,\eta' \in E$, the process
$(e^{-\C t} \cdot \|\hat\Phi_{t,1}(\eta,\eta') - \hat\Phi_{t,2}(\eta,\eta')\|)_{t \ge 0}$
is a supermartingale with respect to its natural filtration.
In particular, for any~$T > 0$ and~$A > 0$, we have
\begin{equation}
\label{eq:first_main_martbound}
\P\left(\sup_{0 \le t \le T}\|\hat\Phi_{t,1}(\eta,\eta')-\hat\Phi_{t,2}(\eta,\eta')\| > A\right) \le \frac{e^{\C T}\cdot \|\eta-\eta'\|}{A}.
\end{equation}
\end{lemma}
\begin{proof}
The first statement is  a consequence of Lemma~\ref{lem:super_from_mc} and Lemma~\ref{lem:bound_couple}. To prove the second statement, abbreviate
\[Y_t:= \|\hat{\Phi}_{t,1}(\eta,\eta') - \hat{\Phi}_{t,2}(\eta,\eta')\|,\qquad X_t := e^{-\C t}\cdot Y_t. \]
For~$a > 0$, define~$\tau_a := \inf\{t \ge 0:\;X_t > a\}$.
We have~$\|\eta-\eta'\| \ge \E[X_{\tau_a \wedge T}] \ge a \cdot \P(\tau_a \le T)$ for any~$T > 0$ and~$a > 0$, so
\[
\P\left(\sup_{0 \le t \le T}X_t > a\right) \le \frac{\|\eta - \eta'\|}{a}.
\]
We then obtain
\[
\P\left(\sup_{0 \le t \le T}Y_t > A \right) \le \P\left(\sup_{0\le t \le T}X_t > Ae^{-\C T}\right) \le \frac{e^{\C T}\cdot \|\eta-\eta'\|}{A}.
\]
\end{proof}

\section{Convergence to solutions of reaction-diffusion equations}
\subsection{Sequence of particle systems: definition and first estimates} \label{ss:sequence}
In this section, following the program outlined in the Introduction, we consider a sequence of processes
of the type constructed in the previous section, and prove that this sequence converges to solutions of the system of reaction-diffusion equations~\eqref{eq:sde}.

We recall that we define, for~$u \ge 0$,
\begin{align*}
  &F^{n,-}(u):= \frac{an^2}{2}\cdot \left(\frac{u}{n}\right)^\ell + \min\left\{\frac{an^2}{2}\cdot
    \left(\frac{u}{n}\right)^\ell ;\;\frac{bn}{2}\cdot\left(\frac{u}{n}\right)^\kappa \right\},\\
  &F^{n,+}(u):= \frac{an^2}{2}\cdot \left(\frac{u}{n}\right)^\ell - \min\left\{\frac{an^2}{2}\cdot
    \left(\frac{u}{n}\right)^\ell ;\;\frac{bn}{2}\cdot\left(\frac{u}{n}\right)^\kappa \right\}.
\end{align*}
We denote by~$(\Phi^n_t(\eta))_{t \ge 0}$ the Markov chains, as constructed in Section~\ref{s:ips},
corresponding to the rate functions~$F^+= F^{n,+}$ and~$F^-=F^{n,-}$
(we leave the diffusion part of the dynamics constant for all~$n$, that is, particles jump with rate one regardless of~$n$).
Conditions in~\eqref{decreasingF} are satisfied with this choice of rate functions.
We also denote by~$\hat{\Phi}^n_\cdot(\eta,\eta')$ the corresponding coupling, as described in the previous section.
Finally, we write
\[
  \varphi^n_t(\zeta) := \tfrac{1}{n}\Phi^n(n\zeta),\quad
  \hat{\varphi}^n_t(\zeta,\zeta')  =\Big(\varphi_{t,1}(\zeta,\zeta'),\hat\varphi_{t,2}(\zeta,\zeta')\Big)
  := \tfrac{1}{n}\hat{\Phi}^n_t(n\zeta,n\zeta'),\quad \zeta,\zeta' \in \tfrac{1}{n}E.
\]
Hence,~$\varphi^n_\cdot(\zeta)$ and~$\hat{\varphi}^n_{\cdot,i}(\zeta,\zeta')$, for $i \in \{1,2\}$, are processes on~$\tfrac{1}{n}E$. 

From~\eqref{eq:first_main_martbound} we obtain that, for any~$n$ and any~$\zeta,\zeta' \in \tfrac{1}{n}E$,
\begin{equation}
\label{eq:second_main_martbound}
\P\left(\sup_{0 \le t \le T}\|\hat\varphi_{t,1}^n(\zeta,\zeta')-\hat\varphi_{t,2}^n(\zeta,\zeta')\| > A\right) \le \frac{e^{\C T}\cdot \|\zeta-\zeta'\|}{A}.
\end{equation}
We have the following consequence of Lemma~\ref{lem:sup0}.
\begin{corollary}\label{lem:eps_delta}
  For every~$\epsilon > 0$ and~$T > 0$, there exists~$\delta > 0$
  such that, for all~$n \in \N$ and all~$\zeta,\zeta' \in \tfrac{1}{n}E$
  with~$\|\zeta-\zeta'\| < \delta$ we have
  \begin{equation}  \label{eq:coupled_zetas}
    \P\left(\sup_{0 \le t \le T}\left\|\hat\varphi^n_{t,1}(\zeta,\zeta') - \hat\varphi^n_{t,2}(\zeta,\zeta')\right\| > \epsilon\right)< \epsilon.
  \end{equation}
\end{corollary}
We may now state a distance bound  with respect to the Lévy-Prohorov metric.
\begin{lemma}\label{cor:uniform}
For every~$\epsilon > 0$, there exists~$\delta > 0$ such that, for all~$n \in \N$
and all~$\zeta,\zeta' \in \tfrac{1}{n}E$ with~$\|\zeta-\zeta'\| < \delta$ we have
\begin{equation}\label{eq:sup_eps_delta}
  \mathrm{d}_{\mathrm{LP}}\left(\varphi^n_\cdot(\zeta),\;\varphi^n_\cdot(\zeta')\right)   < \epsilon.
\end{equation}
\end{lemma}
\begin{proof}
Fix~$\epsilon > 0$. It follows from the definition of the Skorokhod metric that
\[
  \gamma,\gamma' \in D([0,\infty),\mathcal{E}),\; \sup_{0 \le t \le \log(2/\epsilon)}\|\gamma_t - \gamma'_t\|
  \le \frac{\epsilon}{2}
  \quad \Longrightarrow \quad
  \mathrm{d}_\mathrm{S}(\gamma,\gamma')< \epsilon.
\]
Together with Lemma~\ref{lem:eps_delta} (with~$T = \log(2/\epsilon)$),
this implies that there exists~$\delta > 0$ such that for any~$n$ and any~$\zeta,\zeta'$
with~$\|\zeta - \zeta'\| < \delta$ we have
\[
\P\left(\mathrm{d}_\mathrm{S}(\hat{\varphi}^n_{1,\cdot}(\zeta,\zeta'), \hat{\varphi}^n_{2,\cdot}(\zeta,\zeta')) > \epsilon\right) < \epsilon.
\]
The desired result now follows from
\eqref{LPcoup} and the fact that $\hat{\varphi}^n_{1,\cdot}(\eta,\eta')$
and $\hat{\varphi}^n_{2,\cdot}(\eta,\eta')$ have the same distribution of
$\varphi^n_\cdot(\eta)$, $\varphi^n_\cdot(\eta')$ respectively. 
\end{proof}
For a local function~$f:\frac{1}{n}E \to \mathbb{R}$, define
\begin{align*}
&\mathcal{L}^nf(\zeta):= \sum_{\m\in \mathcal{M}} R^\m(n\zeta)\cdot \left(f\left(\tfrac{1}{n}\Gamma^\m(n\zeta)\right) - f(\zeta)\right),\\
&\mathcal{Q}^nf(\zeta):= \sum_{\m\in \mathcal{M}} R^\m(n\zeta)\cdot \left(f\left(\tfrac{1}{n}\Gamma^\m(n\zeta)\right) - f(\zeta)\right)^2,\quad \zeta \in \tfrac{1}{n}E.
\end{align*}
It follows from Lemma~\ref{lem:two_martingales} that the processes
\begin{align*}
  M^f_t & := f(\varphi^n_t(\zeta)) - f(\zeta) - \int_0^t \mathcal{L}^n f(\varphi^n_s(\zeta))ds,
  \quad \text{and}\\
  N^f_t & := \left( M^f_t \right)^2 - \int_0^t \mathcal{Q}^nf(\varphi^n_s(\zeta))ds
\end{align*}
defined for $t \geq 0$ are local martingales.

In the following lemma, we give explicit expressions for~$\mathcal{L}^n$ and~$\mathcal{Q}^n$
applied to functions of the form~$f_x(\zeta) = \zeta(x)$ and~$f_{xy}(\zeta) = \zeta(x)\cdot\zeta(y)$.
We postpone the calculations to Appendix~\ref{sec:app_gen}.
\begin{lemma}\label{lem:expression_generator}
We have, for any~$\zeta \in \frac{1}{n}E$ and~$x \in \V$,
\begin{align}\label{eq:gen_exp_Lnx}
&(\mathcal{L}^nf_x)(\zeta) = (\Delta_p \zeta)(x) - \min\left\{an\cdot (\zeta(x))^\ell; \;b\cdot (\zeta(x))^\kappa\right\}\\
\label{eq:gen_exp_Qnx}
&(\mathcal{Q}^nf_x)(\zeta) = a\cdot (\zeta(x))^\ell + \frac{1}{n}\sum_{y \neq x} (p(y,x)\zeta(y)+p(x,y)\zeta(x)),\\
&\label{eq:gen_exp_Lnxx}(\mathcal{L}^nf_{xx})(\zeta) =  (2f_x \cdot \mathcal{L}^nf_x)(\zeta) + a\cdot (\zeta(x))^\ell + \frac{1}{n}\sum_{y \neq x} (\zeta(y)p(y,x)+\zeta(x)p(x,y)).
\end{align}
Moreover, for distinct~$x,y\in \V$,
\begin{equation}\label{eq:gen_exp_Lnxy}
(\mathcal{L}^nf_{x,y})(\zeta) = (f_x \cdot \mathcal{L}^n f_y + f_y \cdot \mathcal{L}^n f_x)(\zeta) - \frac{1}{n} (\zeta(y) p(y,x)+ \zeta(x)p(x,y)). 
\end{equation}
\end{lemma}

Finally, we state our tightness result, which will allow us to extract convergent subsequences of a sequence of processes of the form~$\{\varphi^n_\cdot(\zeta^n)\}_{n \ge 1}$.
\begin{proposition} \label{prop:tight}
Let~$\zeta  \in \mathcal{E}$ and~$\{\zeta^n\}_{n \ge 1}$
be a sequence with~$\zeta^n \in \tfrac{1}{n}E$ for each~$n$
and $\|\zeta^n - \zeta\| \xrightarrow{n \to \infty} 0$.
Then, the  family of processes~$\{\varphi^n_\cdot(\zeta^n)\}_{n \ge 1}$
is tight in~$D( [0,\infty),\mathcal{E} )$.
\end{proposition}
The proof of this proposition will be carried out in Section~\ref{s:tight}.

\subsection{Limit points are solutions}
Recall from Section~\ref{s:prop_sol} that we defined
\begin{equation*}
(\mathcal{L}^*f)(\zeta) := \sum_{x \in \V} (\Delta_p\zeta(x) - b\cdot (\zeta(x))^\kappa)\cdot \partial_xf(\zeta) + \frac12\sum_{x \in \V} a \cdot (\zeta(x))^\ell\cdot \partial^2_x f(\zeta),\quad \zeta \in \mathcal{E}
\end{equation*}
for any local and twice continuously differentiable function~$f$.
Substituting $f_x$, $f_{xx}$ and~$f_{xy}$ gives
\begin{align}
 \label{eq:gen_exp_L*x}&(\mathcal{L}^*f_x)(\zeta) = \Delta_p\zeta(x) - b\cdot (\zeta(x))^\kappa,\\[.2cm]
&\label{eq:gen_exp_L*xx}(\mathcal{L}^*f_{xx})(\zeta)= (2f_x\cdot \mathcal{L}^*f_x)(\zeta) + a \cdot (\zeta(x))^\ell,\\[.2cm]
&\label{eq:gen_exp_L*xy}(\mathcal{L}^*f_{xy})(\zeta) = (f_x \cdot \mathcal{L}^*f_y + f_y \cdot \mathcal{L}^*f_x)(\zeta)\quad \text{for }x \neq y.
\end{align}

\begin{lemma}\label{lem:uniform_L*}
For any~$f \in \{f_x,f_{xy}:x,y \in \V\}$ and~$A > 0$ we have
\begin{equation}\label{eq:desired_sup_conv}
\sup \left\{|(\mathcal{L}^*f)(\zeta) - (\mathcal{L}^nf)(\zeta)| :\; \zeta \in \tfrac{1}{n}E,\; \|\zeta\| \le A   \right\} \xrightarrow{n \to \infty}0.
\end{equation}
\end{lemma}
\begin{proof}
Let us fix $x \in \V$, and consider the case~$f = f_x$.
Comparing~\eqref{eq:gen_exp_Lnx} and~\eqref{eq:gen_exp_L*x},
and noting that~$\zeta(x) \le \|\zeta\|/\alpha(x)$,
the supremum on the left-hand side of~\eqref{eq:desired_sup_conv} is at most
\begin{align*}
& \sup_{z \in [0,A/\alpha(x)]} \left(b z^\kappa - an z^\ell\right) \cdot \mathds{1}\left\{b z^k > an  z^\ell\right\} \xrightarrow{n \to \infty}0;
\end{align*}
the convergence can be checked by separately considering the cases~$\kappa > \ell$ and~$\kappa \le \ell$.

Next, for~$f = f_{xx}$, comparing~\eqref{eq:gen_exp_Lnxx} and~\eqref{eq:gen_exp_L*xx}
and using the case of~$f_x$ that we just treated, it suffices to note that, by~\eqref{eq:quick_comp_q}, we have 
\[
  \frac{1}{n}\cdot \sup \left\{\sum_{y \neq x}(\zeta(y)p(y,x) + \zeta(x)p(x,y)): \;\zeta \in \tfrac{1}{n}E,\;\|\zeta\| \le A \right\} \xrightarrow{n \to \infty}0.
\]

Finally, the case~$f=f_{xy}$ with given~$x\neq y$ is easier:
comparing~\eqref{eq:gen_exp_Lnxy} and~\eqref{eq:gen_exp_L*xy}
and using the convergence for~$f_x$ and~$f_y$, it suffices to note that
\[
  \frac{1}{n}\cdot \sup  \left\{\zeta(y)p(y,x) + \zeta(x)p(x,y):\; \zeta \in \tfrac{1}{n}E,\;\|\zeta\| \le A  \right\}  \xrightarrow{n \to \infty} 0.
  \qedhere
\]
\end{proof}

\begin{proposition}
\label{prop:solutions}
Let~$\zeta \in \mathcal{E}$ and~$\{\zeta^n\}_{n \ge 1}$ be a sequence with~$\zeta^n \in \tfrac{1}{n}E$ for each~$n$ and~$\|\zeta^n - \zeta\| \xrightarrow{n \to \infty}0$. Assume that the sequence of processes~$\{\varphi^{n}_\cdot(\zeta^n)\}_{n \ge 1}$ has a subsequence that converges in distribution in~$D([0,\infty),\mathcal{E})$ to a process~$\zeta_\cdot^*$. Then, the distribution of~$\zeta_\cdot^*$ is supported on~$C([0,\infty),\mathcal{E})$, and~$\zeta_\cdot^*$ is a solution of the SDE~\eqref{eq:sde}. 
\end{proposition}
\begin{proof}[Proof of Proposition~\ref{prop:solutions}] 
Denote the convergent subsequence by~$\{\varphi^{n_k}_\cdot(\zeta^{n_k})\}_{k \ge 1}$.
By Skorokhod's representation theorem~\cite[p.~70]{Bil99}, we may consider a probability space~$(\widetilde{\Omega},\widetilde{\mathcal{F}},\widetilde{P})$ with processes~$\{Z^k_\cdot\}_{k \ge 1}$ and~$Z^*_\cdot$ so that
\begin{itemize}
\item for each~$k$,~$Z^k_\cdot$ has trajectories in~$D([0,\infty),\frac{1}{n}E)$ and same distribution as~$\varphi^{n_k}_\cdot(\zeta^{n_k})$;
\item $Z^*_\cdot$ has trajectories in~$D([0,\infty),\mathcal{E})$ and same distribution as~$\zeta^*_\cdot$;
\item for each~$\omega \in \widetilde{\Omega}$, we have~$\mathrm{d}_\mathrm{S}(Z^k_\cdot(\omega),Z^*_\cdot(\omega))\xrightarrow{n \to \infty}0.$
\end{itemize}

For~$\gamma \in D([0,\infty),\mathcal{E})$ and~$t \ge 0$, define~$J_t(\gamma):= \sup_{s \le t}\|\gamma_s - \gamma_{s-}\|$, the largest jump size of~$\gamma$ until time~$t$.  We have that
\[
  J_t\left(Z^k_\cdot(\omega)\right) \le \frac{2}{n_k}\cdot \max_{x \in \V}\alpha(x) \xrightarrow{k \to \infty} 0\quad \text{ for all }\omega \in \widetilde{\Omega},
\]
so, by continuity of~$J_t$ in the Skorokhod topology~\cite[p.~125]{Bil99},
we obtain~$J_t\left(Z^*_t(\omega)\right) = 0$ for all~$\omega$ and~$t$.
This implies that the trajectories of~$Z^*_\cdot$ are continuous,
and we proved the first part.

Then, using~\eqref{eq:ds_cont}, we obtain
\begin{equation}\label{eq:zs_are_near}
\sup_{0\le s \le t} \|Z^k_s(\omega)-Z^*_s(\omega)\| \xrightarrow{k \to \infty} 0\quad \text{for all }t \ge 0,\;\omega \in \widetilde{\Omega}.
\end{equation}

We now define, for~$f \in \{f_x,f_{xy}:x,y \in \V\}$,
\begin{align*}
&M^{k,f}_t := f(Z^k_t) - f(Z^k_0) - \int_0^t \mathcal{L}^{n_k}f(Z^k_s)ds,\\
&M^{*,f}_t := f(Z^*_t) - f(Z^*_0) - \int_0^t \mathcal{L}^*f(Z^*_s)ds,\quad t \ge 0.
\end{align*}
As observed in Section~\ref{ss:sequence}, we have that~$M^{k,f}_\cdot$ is a local martingale for each~$k$. We now claim that
\begin{equation}\label{eq:sup_of_ms}
\sup_{0\le s \le t} \left| M^{k,f}_s(\omega) - M^{*,f}_s(\omega)\right| \xrightarrow{k \to \infty} 0 \quad \text{for all } t \ge 0,\; \omega \in \widetilde{\Omega}.
\end{equation}
Let us first show how this convergence will allow us to conclude. Since the trajectories of~$M^{*,f}_\cdot$ are continuous,~\eqref{eq:ds_cont} and~\eqref{eq:sup_of_ms} imply that
\[
  \mathrm{d}_\mathrm{s}(M_\cdot^{k,f}(\omega),M_\cdot^{*,f}(\omega))\xrightarrow{k \to \infty} 0 \quad \text{for all~$\omega \in \widetilde{\Omega}$}.
\]
Since almost sure convergence implies convergence in distribution, this  gives
\[
  \mathrm{d}_{\mathrm{LP}}(M^{k,f}_\cdot,M^{*,f}_\cdot)\xrightarrow{k \to \infty} 0.
\]
Now, Corollary~1.19, page 527 in~\cite{shyriaev} states that if a sequence of
c\`adl\`ag local martingales converges in distribution (with respect to the Skorokhod topology),
then the limiting process is also a local martingale.
We then obtain that~$M^{*,f}_\cdot$ is a local martingale.
By Proposition~\ref{prop:are_sols}, this implies  that~$Z^*_\cdot$ is a solution to the SDE~\eqref{eq:sde}.

It remains to prove~\eqref{eq:sup_of_ms}.
To do so, fix~$f \in \{f_x,f_{xy}:x,y \in \mathbb{V}\}$,
$\omega \in \widetilde{\Omega}$ and~$t \ge 0$.
Using~\eqref{eq:zs_are_near}, we obtain that
\[
  \sup_{0 \le s \le t} \Big| f(Z^k_s(\omega)) - f(Z^*_s(\omega)) \Big|
    \xrightarrow{n \to \infty}0,
\]
and also that~$A:= \sup\{\|Z^k_s(\omega)\|: s \le t,\;k \ge 1\}$ is finite.
By this latter point and Lemma~\ref{lem:uniform_L*}, we then have that
\[
  \sup_{0\le s \le t}|\mathcal{L}^{n_k}f(Z^k_s(\omega)) - \mathcal{L}^*f(Z^k_s(\omega))|\xrightarrow{k \to \infty}0.
\]
Next, from the generator expressions in~\eqref{eq:gen_exp_L*x},~\eqref{eq:gen_exp_L*xx}
and~\eqref{eq:gen_exp_L*xy}, it follows that~$\mathcal{L}^*f$
is uniformly continuous on~$\{\zeta \in \mathcal{E}: \|\zeta\| \le A\}$;
this and~\eqref{eq:zs_are_near} imply that
\[
  \sup_{0 \le s \le t} \Big| \mathcal{L}^*f(Z^k_s(\omega)) - \mathcal{L}^*f(Z^*_s(\omega)) \Big|
    \xrightarrow{k \to \infty} 0.
\]
The desired  convergence~\eqref{eq:sup_of_ms} follows.
\end{proof}

\subsection{Convergence to solutions: proof of Theorem~\ref{thm:main}}

\begin{proof}[Proof of Theorem~\ref{thm:main}]
We split the proof in two parts: first with finite initial condition,
then the general case.

\def\ozeta{\mathring{\zeta}}
\paragraph{Finite case}
Fix~$\ozeta \in \mathcal{E}$ with~$\|\ozeta\|_1 = \sum_{x\in \V} \ozeta(x) < \infty$.
Also, fix a sequence~$\{\ozeta^n\}_{n \ge 1}$ with
$\ozeta^n \in \tfrac{1}{n}E$ for each~$n$
and~$\|\ozeta^n - \ozeta\| \xrightarrow{n \to \infty} 0$.
By Proposition~\ref{prop:tight},
there exists a subsequence~$\{\ozeta^{n_k}\}_{k \ge 1}$ such that
the sequence of processes~$\{\varphi^{n_k}_\cdot(\ozeta^{n_k})\}_{k \ge 1}$
converges in distribution to a process on~$D([0,\infty),\mathcal{E})$.
Let us denote this limiting process by~$\psi^*_\cdot(\ozeta)$.
By Proposition~\ref{prop:solutions},
$\psi^*_\cdot(\ozeta)$ has trajectories in~$C([0,\infty),\mathcal{E})$,
and is a solution to~\eqref{eq:sde}
with initial configuration~$\bar{\zeta} = \ozeta$.
Then, by Proposition~\ref{prop:unique},
any other subsequence of~$\{\varphi^n(\ozeta^n)\}_{n \ge 1}$
which converges in distribution must have the same limit~$\psi^*_\cdot(\ozeta)$.
This implies that the entire sequence~$\{\varphi^n(\mathring{\eta}^n)\}_{n \ge 1}$
converges in distribution to~$\psi^*_\cdot(\ozeta)$, that is,
\begin{equation*}
\mathrm{d}_{\mathrm{LP}}\left( \varphi^n_\cdot(\ozeta^n), \;\psi^*_\cdot(\ozeta)\right)\xrightarrow{n \to \infty} 0.
\end{equation*}

\paragraph{General case}
Now, fix~${\zeta} \in \mathcal{E}$ and a sequence~$\{\zeta^n\}_{n \ge 1}$
with~$\zeta^n \in \tfrac{1}{n}E$ for each~$n$
and~$\|{\zeta}^n - {\zeta}\| \xrightarrow{n \to \infty} 0$.
We claim that the sequence of processes~$\{\varphi^n_\cdot(\zeta^n)\}_{n \ge 1}$
is Cauchy with respect to~$\mathrm{d}_\mathrm{LP}$.
To see this, fix~$\epsilon > 0$.
Next, choose~$\delta > 0$ such that Lemma~\ref{lem:eps_delta} ensures that
the left-hand side in~\eqref{eq:sup_eps_delta}
is smaller than~$\frac{\epsilon}{3}$.
Next, because~$\|\zeta^n - \zeta\| \to 0$, we may choose~$R > 0$ such that
\[
  \|\zeta^n \cdot \mathds{1}_{B_0(R)} - \zeta^n\| < \delta\quad \text{ for all~$n$}.
\]
This implies that
\[
  \mathrm{d}_{\mathrm{LP}} \left(\varphi^n_\cdot(\zeta^n \cdot \mathds{1}_{B_0(R)}),\;
    \varphi^{n}_\cdot(\zeta^{n})\right) < \frac{\epsilon}{3}
  \quad \text{for all $n \in \N$}.
\]
Finally, using the finite case, with $\{\zeta^n \cdot \mathds{1}_{B_0(R)}\}_{n}$
as approximating sequence to $\zeta \cdot \mathds{1}_{B_0(R)}$, which is finite,
we may choose $n_0 \in \N$ such that
\[
  \mathrm{d}_{\mathrm{LP}}\left(\varphi^n_\cdot(\zeta^n \cdot \mathds{1}_{B_0(R)}),\;
    \varphi^{n'}_\cdot(\zeta^{n'} \cdot \mathds{1}_{B_0(R)})\right)< \frac{\epsilon}{3}\quad \text{for all }n,n' \ge n_0.
\]
By the triangle inequality, it follows that~$\{\varphi^n_\cdot(\zeta^n)\}_{n \ge 1}$ is Cauchy. 

Now, from the tightness of this sequence,
given by Proposition~\ref{prop:tight}
(or alternatively, the completeness of
the metric space~$(D([0,\infty),\mathcal{E}),\mathrm{d}_{\mathrm{LP}}) $),
we obtain that~$\{\varphi^n_\cdot(\eta^n)\}_{n \ge 1}$ converges in distribution
to a process on~$D([0,\infty),\mathcal{E})$,
which we denote by~$\psi^*_\cdot(\zeta)$.
To see that this process does not depend
on the sequence~$\{\zeta^n\}_{n \ge 1}$ that we fixed,
take an alternative sequence~$\{\tilde{\zeta}^n\}_{n \ge 1}$
with~$\tilde{\zeta}^n \in \frac{1}{n}E$ for each~$n$
and~$\|\tilde{\zeta}^n - \zeta\| \xrightarrow{n \to \infty}0$,
and note that Lemma~\ref{lem:eps_delta} gives
\[
  \lim_{n \to \infty} \mathrm{d}_{\mathrm{LP}}\left(\varphi^n_\cdot(\zeta^n),\;
    \varphi^n_\cdot(\tilde{\zeta}^n) \right) \xrightarrow{n \to \infty}0.
\]
Finally, by Proposition~\ref{prop:solutions},
$\psi^*_\cdot({\zeta})$ has trajectories in~$C([0,\infty),\mathcal{E})$,
and is a solution to~\eqref{eq:sde}
with initial configuration~$\bar{\zeta} = {\zeta}$.

It remains to prove part~(c) of the statement of the theorem.
Fix~$\epsilon > 0$, and choose~$\delta > 0$
corresponding to~$\epsilon$ in Lemma~\ref{cor:uniform}.
Fix~$\zeta,\zeta' \in \mathcal{E}$ with~$\|\zeta - \zeta'\| < \delta/2$.
Take sequences~$\{\zeta^n\}_{n\ge 1}$, $\{\zeta'^n\}_{n \ge 1}$
with~$\zeta^n,\zeta'^n \in \tfrac{1}{n}E$ for each~$n$ and
\[
  \|\zeta^n - \zeta\| \xrightarrow{n \to \infty} 0,
  \qquad
  \|\zeta'^n - \zeta'\| \xrightarrow{n \to \infty} 0.
\]
In particular, for~$n$ large enough we have~$\|\zeta^n - \zeta'^n\| < \delta$,
so $\mathrm{d}_{\mathrm{LP}}(\varphi^n_\cdot(\zeta),\varphi^n_\cdot(\zeta')) < \epsilon$.
By the previous results, we have
\[
  \mathrm{d}_{\mathrm{LP}}(\varphi^n_\cdot(\zeta^n),\psi^*_\cdot(\zeta)) \xrightarrow{n \to \infty} 0,
  \qquad
  \mathrm{d}_{\mathrm{LP}}(\varphi^n_\cdot(\zeta'^n),\psi^*_\cdot(\zeta')) \xrightarrow{n \to \infty} 0,
\]
so
\[
  \mathrm{d}_{\mathrm{LP}}(\psi^*_\cdot(\zeta),\psi^*_\cdot(\zeta'))
  = \lim_{n \to \infty} \mathrm{d}_{\mathrm{LP}}(\varphi^n_\cdot(\zeta^n),\varphi^n_\cdot(\zeta'^n))
  \le \epsilon.
  \qedhere
\]
\end{proof}

\section{Tightness: proof  of Proposition~\ref{prop:tight}} \label{s:tight}
\subsection{Aldous' criterion}
Throughout this section, we fix~$\zeta$ and~$\{\zeta^n\}_{n \ge 1}$
as in the statement of Proposition~\ref{prop:tight}. We will abbreviate
\[
  \zeta^n_t:= \varphi^n_t(\zeta^n),\quad n \ge 1,\; t \ge 0.
\]

We will prove Proposition~\ref{prop:tight} using Aldous' criteria~\cite[p. 51]{Bil99}. We need to verify that
\begin{align}
  \label{eq:aldous1}
  &\forall t \ge 0, \;\forall \epsilon > 0,\; \exists K \subset \mathcal{E} \text{ compact}: \;\sup_{n \in \N}\;
    \P(\zeta^n_t \notin K) < \epsilon, \text{ and}\\
  \label{eq:aldous2}
  &\forall T > 0,\;\forall \epsilon > 0,\quad \lim_{\delta \to 0}\;\sup_{n \in \N}\;\sup_{\tau \in \mathcal{T}^n_T}\;
    \P\left(\left\|\zeta^n_{(\tau+\delta)\wedge T} - \zeta^n_\tau \right\|  > \epsilon \right) = 0,
\end{align}
where~$\mathcal{T}^n_T$ is the set of stopping times (with respect to the natural filtration of~$\zeta^n_\cdot$) that are bounded by~$T$. 

To verify the first criterion, we will rely on a definition and a lemma. For any~$r > 0$ we define
\begin{equation}\label{eq:def_Lambda}\Lambda(r):= \{x \in \V: \;\alpha(x) > 1/r\}. \end{equation}
\begin{lemma}[Negligible norm near infinity]\label{lem:neg}
For any~$T > 0$ and~$\epsilon > 0$ there exists~$R > 0$ such that
\begin{equation}
\P\left( \sup_{0 \le t \le T} \|\zeta^n_t \cdot \mathds{1}_{\Lambda(R)^c}\| > \epsilon\right) < \epsilon \quad \text{for all }n \ge 1.
\end{equation}
\end{lemma}
We postpone the proof of this  lemma to Section~\ref{ss:neg}.
\begin{proof}[Proof of Proposition~\ref{prop:tight}, condition~\eqref{eq:aldous1}]
Fix~$t \ge 0$ and~$\epsilon > 0$. Using~\eqref{eq:second_main_martbound} with~$\zeta' \equiv 0$,
we obtain that there exists~$A > 0$ such that
\[
  \sup_n \P\left(\|\zeta^n_t\| > A\right) \le \frac{\epsilon}{2}.
\]
Furthermore, by Lemma~\ref{lem:neg}, for any~$k \in \N$ there exists~$R_k$ such that
\[
  \sup_n \P\left(\|\zeta^n_t \cdot \mathds{1}_{\Lambda(R_k)^c} \| > \frac{1}{k} \right)< \frac{\epsilon}{2^{k+1}}.
\]
Now, defining
\[
  K:= \left\{\zeta \in \mathcal{E}:\;\|\zeta\| \le A,\; \|\zeta\cdot \mathds{1}_{\Lambda(R_k)^c}\| \le 1/k \text{ for all }k  \right\},
\]
we have that~$\P(\zeta_t^n \in K) > 1-\epsilon$ for all~$n$. 

We claim that~$K$ is compact. To verify this, fix a sequence~$\{\zeta^j\}_{j \ge 1}$ of elements of~$K$. For every~$x \in \V$ we have that~$\zeta^j(x) \le A/\alpha(x)$ for every~$j$, so, using a diagonal argument, we can obtain a subsequence~$\{\zeta^{j'}\}_{j' \ge 1}$ so that~$\zeta^{j'}(x)$ is convergent for each~$x$. Let~$\bar{\zeta}$ be defined by~$\bar{\zeta}(x):= \lim_{j'} \zeta^{j'}(x)$ for each~$x$. Next, Fatou's Lemma gives~$\|\bar\zeta \cdot \mathds{1}_{\Lambda(R_k)^c}\| \le 1/k$ for all~$k$, so~$\bar\zeta \in \mathcal{E}$ (since~$\Lambda(R_k)$ is finite) and
\[
  \limsup_{j' \to \infty} \|\bar\zeta - \zeta^{j'}\| \le \limsup_{j' \to \infty} \sum_{x \in \Lambda(R_k)}\alpha(x)\cdot |\bar \zeta(x) - \zeta^{j'}(x)| + \frac2k = \frac2k,
\]
which can be made as small as desired by  taking~$k$ large. This shows that~$K$ is compact and completes  the proof of~\eqref{eq:aldous1}.
\end{proof}

For the proof of~\eqref{eq:aldous2}, again we will need a preliminary result.
\begin{lemma}[Oscillation of coordinates] \label{lem:osc}
For any~$T > 0$,~$\epsilon > 0$ and~$x \in \V$, we have
\begin{equation}
\lim_{\delta \to 0}\;\sup_{n \in \N}\;\sup_{\tau  \in \mathcal{T}_T^n} \; \P\left( \left| \zeta^n_{(\tau+\delta)\wedge T}(x) - \zeta^n_\tau(x)  \right| > \epsilon\right) = 0.
\end{equation}
\end{lemma}
We postpone the proof of this lemma to Section~\ref{ss:osc}.
\begin{proof}[Proof of Proposition~\ref{prop:tight}, condition~\eqref{eq:aldous2}] Fix~$T  > 0$ and~$\epsilon > 0$. By Lemma~\ref{lem:neg}, we may choose~$R> 0$ such that, for any~$n$,
  \[
    \P\left(\sup_{0\le t \le T} \| \zeta^n_t \cdot \mathds{1}_{\Lambda(R)^c}\| > \frac{\epsilon}{3}\right)<\frac{\epsilon}{3}.
  \]
Next, by Lemma~\ref{lem:osc} and a union bound for~$x \in \Lambda(R)$, we can choose~$\delta_0 > 0$ such that, for any~$\delta \le \delta_0$,~$n \in \N$ and~$\tau \in \mathcal{T}_T^n$, we have
\[
  \mathbb{P}\left( \| (\zeta^n_{(\tau+\delta)\wedge T} - \zeta^n_t)\cdot \mathds{1}_{\Lambda(R)} \| > \frac{\epsilon}{3} \right)<\frac{\epsilon}{3}.
\]
Therefore, by the triangle inequality,
\begin{align*}
\P\left( \| \zeta^n_{(\tau + \delta)\wedge T} - \zeta^n_t \| > \epsilon \right) \le & \P\left(\|(\zeta^n_{(\tau+\delta)\wedge T}- \zeta^n_\tau) \cdot \mathds{1}_{\Lambda(R)}\| > \frac{\epsilon}{3} \right) \\&+\P\left( \|\zeta^n_{(\tau+\delta)\wedge T} \cdot \mathds{1}_{\Lambda(R)^c}  \| > \frac{\epsilon}{3}\right) + \P\left( \| \zeta^n_\tau \cdot \mathds{1}_{\Lambda(R)^c}\| > \frac{\epsilon}{3} \right)< \epsilon.
\end{align*}
Since~$\epsilon$ is arbitrary, the proof is complete.
\end{proof}

\subsection{Norm near infinity: proof of Lemma~\ref{lem:neg}}\label{ss:neg}
Let us define
\[
  \zeta^{n,r}_t:= \varphi^n_t(\zeta^n\cdot \mathds{1}_{\Lambda(r)}),\quad n \ge 1,\;t \ge 0.
\]
We observe that, by~\eqref{eq:second_main_martbound}, for any~$\epsilon > 0$,~$T > 0$ and~$n \ge 1$ we have
\begin{equation}\label{eq:apply_unif}
 \P\left( \sup_{0 \le t \le T} \|\zeta^n_t - \zeta^{n,r}_t\| > \epsilon\right) \le \frac{e^{\C T}\cdot \| \zeta_0^n \cdot \mathds{1}_{\Lambda(r)^c}\|}{\epsilon}.
\end{equation}
\begin{proof}[Proof of Lemma~\ref{lem:neg}]
Fix~$T > 0$ and~$\epsilon > 0$. By~\eqref{eq:apply_unif}, we can choose~$r$ large enough that
\begin{equation}
\label{eq:apply_unif2}
\sup_n \P\left(\sup_{0\le t \le T} \|\zeta^n_t - \zeta^{n,r}_t\| > \frac{\epsilon}{2}\right) \le \frac{\epsilon}{2}.
\end{equation}
Next, note that for any~$n$,~$R$ and~$t$ we have
\[
  \|\zeta^{n,r}_t \cdot \mathds{1}_{\Lambda(R)^c}\| = \sum_{x \notin \Lambda(R)} \alpha(x)\cdot |\zeta^{n,r}_t(x)|
  \stackrel{\eqref{eq:def_Lambda}}{\le} \frac{\|\zeta^{n,r}_t\|_1}{R},
\]
so, for any~$n$ and~$R$,
\begin{equation}\label{eq:apply_unif3}
  \P\left(\sup_{0 \le t \le T} \|\zeta^{n,r}_t \cdot \mathds{1}_{\Lambda(R)^c} \| > \frac{\epsilon}{2} \right)
  \le \P\left( \sup_{0\le t \le T} \|\zeta^{n,r}_t\|_1 > \frac{R\epsilon}{2} \right)
  \stackrel{\eqref{eq:l1norm}}{\le}\frac{2\|\zeta^{n,r}_0\|_1}{R\epsilon}.
\end{equation}
Now, the assumption that~$\|\zeta_0^n -\zeta^*\| \xrightarrow{n \to \infty}0$ implies that~$\sup_n \|\zeta_0^{n,r}\|_1 < \infty$. We thus have
\begin{equation}\label{eq:apply_unif4}
\sup_n \P \left(\sup_{0 \le t \le T}\|\zeta^{n,r}_t \cdot \mathds{1}_{\Lambda(R)^c}\|  > \frac{\epsilon}{2} \right) \le \frac{\epsilon}{2}
\end{equation}
if~$R$ is large enough. Combining~\eqref{eq:apply_unif2} and~\eqref{eq:apply_unif4} with the bound
\[
  \|\zeta^n_t \cdot \mathds{1}_{\Lambda(R)^c} \|
  \le \|(\zeta^n_t - \zeta^{n,r}_t)\cdot \mathds{1}_{\Lambda(R)^c} \| + \|\zeta^{n,r}_t\cdot \mathds{1}_{\Lambda(R)^c} \|
  \le \|\zeta^n_t - \zeta^{n,r}_t\| + \|\zeta^{n,r}_t\cdot \mathds{1}_{\Lambda(R)^c} \|
\]
gives the desired bound.
\end{proof}
\subsection{Oscillation of coordinates: proof of Lemma~\ref{lem:osc}}
\label{ss:osc}
In the proof of Lemma~\ref{lem:osc}, it will be useful to note that, for any~$n \in \N$,~$x \in \V$ and~$A > 0$ we have
\begin{equation}\label{eq:sup_of_max}
  \sup\left\{|\mathcal{L}^nf_x(\zeta)|\vee |\mathcal{Q}^nf_x(\zeta)|:\;
    \zeta \in \tfrac{1}{n}E,\; \|\zeta\| \le A\right\} < \infty.
\end{equation}
This follows from the expressions in~\eqref{eq:gen_exp_Lnx}, \eqref{eq:gen_exp_Qnx},
the fact that~$\|\zeta\| \le A$ implies~$|\zeta(x)| \le A / \alpha(x)$,
and the bound
\begin{equation}
  \label{eq:quick_comp_q}
  \sum_{y \neq x}(\zeta(y)p(y,x) + \zeta(x)p(x,y))
  \stackrel{\eqref{eq:quick_bound}}{\le}
  \sum_{y \in \V} \zeta(y) \cdot \frac{\C\alpha(y)}{\alpha(x)} + \zeta(x)
  \le \frac{\C+1}{\alpha(x)}\cdot \|\zeta\|.
\end{equation}
\begin{proof}[Proof of Lemma~\ref{lem:osc}]
As noted in Section \ref{ss:sequence}, writing
\begin{equation}
  \label{e:cm}
  \mc{M}^{n,x}_t :=  \X^n_t(x) - \X^n_0(x) - \int_0^{t}\mc{L}^nf_{x}(\X^n_s)\, ds,
  \qquad \text{for $t \ge 0$},
\end{equation}
we have that~$\mc{M}^{n}_\cdot$ is a local martingale.
Its quadratic variation is given by
\begin{equation}\label{e:cqv}
  \<\mc{M}^{n,x}\>_t
  = \int_0^{t}\mc{Q}^nf_{x}(\X^n_s)\, ds
  = \int_0^t a \cdot (\X^n_s(x))^\ell + \sum_{y} \frac{p(y,x)\X^n_s(y) + p(x,y)\X^n_s(x)}{n} \, ds,
  \quad \text{for $t \ge 0$}.
\end{equation}

Now, fix~$T > 0$, $\epsilon > 0$ and~$\tau \in \mathcal{T}^n_T$.
Given~$A > 0$, define $\tau^n_A := T \wedge \inf \chv{t \in [0,T]: \norm{\X^n_t} > A}$.
We have that, for any~$\delta > 0$,
\begin{equation}
\label{oscilado}
\begin{aligned}
&\bb{P}\prt{\abs{\X^n_{(\tau+\delta)\wedge T}(x) - \X^n_\tau(x)} \geq \gep}\\
&\leq \bb{P}\prt{\abs{\mc{M}^n_{(\tau+\delta)\wedge \tau^n_A} - \mc{M}^n_{\tau \wedge \tau^n_A}} \geq \gep/2}
+ \bb{P}\prt{\abs{\int_{\tau\wedge \tau^n_A}^{(\tau + \delta) \wedge \tau^n_A}\mc{L}^nf_{x}(\X^n_s)\,ds} \geq \gep/2} + \bb{P}(\tau^n_A <T)\\
&\leq  \frac{4}{\gep^2} \bb{E}\crt{\int_{\tau\wedge \tau^n_A}^{(\tau + \delta)\wedge \tau^n_A}\mc{Q}^nf_x(\X^n_s)\, ds}
+ \bb{P}\prt{\abs{\int_{\tau\wedge \tau^n_A}^{(\tau + \delta) \wedge \tau^n_A}\mc{L}^nf_x(\X^n_s)\,ds} \geq \gep/2}
 + \bb{P}(\tau^n_A <T).
\end{aligned}
\end{equation}
Now, using~\eqref{eq:sup_of_max}, we can choose~$\delta > 0$ small enough that
\[
  \frac{4}{\gep^2} \bb{E}\crt{\int_{\tau\wedge \tau^n_A}^{(\tau + \delta)\wedge \tau^n_A}\mc{Q}^nf_x(\X^n_s)\, ds} < \frac{\epsilon}{2}
  \qquad  \text{and}\qquad
  \bb{P}\prt{\abs{\int_{\tau\wedge \tau^n_A}^{(\tau + \delta) \wedge \tau^n_A}\mc{L}^nf_x(\X^n_s)\,ds} \geq \gep/2} = 0.
\]
To conclude the proof, we choose~$A > 0$ large enough
that~$\P(\tau^n_A < T) < \tfrac{\epsilon}{2}$ .
\end{proof}

\section*{Acknowledgements}
Motivated by a question posed by Roberto Oliveira, this paper started during the first author postdoc at IMPA (BR).
Collaboration initiated with visits at the Mathematical Institute in Leiden and at the Bernoulli Institute in Groningen University (NL).
Important progress happened while the second author visited the Department of Mathematical Sciences in Durham University (UK),
while on leave from Federal University of Rio de Janeiro (BR).
We would like to express our gratitude for the hospitality and support of all those institutions.

\appendix

\section{Appendix}\label{s:appendix}
\subsection{Uniqueness with finite mass: proof of Proposition~\ref{prop:unique}} \label{s:unique}

Let $\X^{*,1}_\cdot$ and $\X^{*,2}_\cdot$ be two limit points of $\{\X^n_\cdot\}_{n \in \V}$
having the same initial condition $\X^*_0 \in \mc{E}_0$,
By~\cite[Proposition 1, p. 158]{WatYam71}, we can assume that the two processes~$\X^{1}_\cdot$ and $\X^{2}_\cdot$
of the statement of the proposition are defined in the same probability space,
and that they solve~\eqref{eq:sde} with respect to the same Brownian family $\chv{B^x_\cdot}_{x \in \V}$.
This allows us to consider the difference process $D_t(x) := \X^{*,1}_t(x) - \X^{*,2}_t(x)$.
To prove that $\X^{1}_\cdot$ and $\X^{2}_\cdot$ have the same distribution
it is enough to show that, for all $t \in [0,T]$,
\begin{equation}\label{equalmass1}
\bb{E}\crt{\norm{D_t}_1} = 0,
\end{equation}
where $\norm{D_t}_1 = \sum_x \abs{D_t(x)}$.
To do so, we first bound $\bb{E}\crt{\norm{D_{t\wedge \sigma_A}}_1}$,
where
\begin{equation}\label{sigmastop}
\sigma_A := \inf\set{t : \|\X^{*,1}_t\|_1 \geq A \text{ or } \|\X^{*,2}_t\|_1 \geq A}.
\end{equation}

As in~\cite[p.165]{WatYam71}, there exist a decreasing sequence $(a_m)_{m \ge 1} $
of positive real numbers with $\lim_m a_m = 0$ and functions~$\rho_m: [0,\infty) \to [0,\infty)$ such that
\begin{equation}\label{rhocond}
  0 \leq \rho_m(u) \leq \frac{2}{mu}\; \forall u,
  \quad
  \rho_m(x) = 0 \text{ for } x \notin (a_m, a_{m-1}),
  \quad
  \text{ and }
  \quad
  \int_{a_m}^{a_{m-1}} \rho_m(s)\, ds = 1.
\end{equation}
We then define
\begin{equation}\label{c2phi}
\varphi_m(u) := \int_0^{\abs{u}}\int_0^y \rho_m(s) \, ds dy, \qquad u \in \R.
\end{equation}
Note that~$\varphi_m$ is twice differentiable, and moreover we have
$\abs{\varphi'_m(u)}\leq 1$, $\varphi_m \leq \varphi_{m+1}$, and~$\lim_m \varphi_m(u) = \abs{u}$.

By Itô's formula,
\begin{equation}\label{e:itouniq}
\begin{aligned}
\varphi_m(D_{t\wedge \sigma_A}(x))
 = &\int_0^{t \wedge \sigma_A} \varphi_m'(D_s(x)) \crt{\Delta_p D_s(x) -b \cdot \prt{\prt{\X^{1}_s(x)}^\kappa - \prt{\X^{2}_s(x)}^\kappa}}\, ds\\
& + \frac{1}{2} \int_0^{t\wedge \sigma_A} \varphi_m''(D_s(x)) \prt{a \cdot \prt{\prt{\X^{1}_s(x)}^{\ell/2} - \prt{\X^{2}_s(x)}^{\ell/2}}}^2\, ds\\
& + \int_0^{t\wedge \sigma_A} \varphi_m'(D_s(x)) \prt{a \cdot \prt{\prt{\X^{1}_s(x)}^{\ell/2} -\prt{\X^{2}_s(x)}^{\ell/2}}}\, dB^x_s.
\end{aligned}
\end{equation}
The expected value of the last integral is zero.
To bound the first two, we note that for every~$s \leq \sigma_A$ and~$\i \in \V$
we have~$\X^{1}_s(\i) \leq A$ and~$\X^{2}_s(\i) \leq A$. In the remainder of this proof,~$C_A$
will denote a positive constant that depends only on~$A$, and whose value may change from line to line.
Since both $\kappa$ and $\ell$ are larger than or equal to one, we have 
\begin{align}
  \label{kdifbound}
  \prt{\X^{1}_s(x)}^\kappa - \prt{\X^{2}_s(x)}^\kappa & \leq C_A \abs{D_s(x)}, \text{ and}\\
  \label{elldifbound}
  \prt{\prt{\X^{1}_s(x)}^{\ell/2} - \prt{\X^{2}_s(x)}^{\ell/2}}^2 & \leq C_A \abs{D_s(x)}.
\end{align}
Since
\begin{equation}\label{absbondrho}
\varphi''_m(D_s(x)) = \rho_m(\abs{D_s(x)}) \leq \frac{2}{m\abs{D_s(x)}},
\end{equation}
the second integral in \eqref{e:itouniq} is bounded by $\frac{C_A a t}{m}$,
so
\begin{equation}
  \label{e:lipbound}
  \bb{E}\crt{\varphi_m(D_{t\wedge \sigma_A}(x))}
  \leq \bb{E}\crt{\int_0^{t\wedge \sigma_A} \abs{\varphi_m'(D_s(x))}
    \crt{\strut\abs{\Delta_p D_s(x)} + C_A \abs{D_s(x)}}\, ds} + \frac{C_A a t}{m}.
\end{equation}

Now define $v_t(\i) := \bb{E}\crt{\abs{D_{t \wedge \sigma_A}(\i)}}$ for all $x \in V$
and let $m \to \infty$ in~\eqref{e:lipbound}.
Since $\sup_u\abs{\varphi'_m(u)} \leq 1$, it follows that
\begin{equation}\label{siteEbound}
v_t(\i) \leq \int_0^{t} \Big(\sum_{\j\in \V } p(\j,\i) v_s (\j) + p(\i,\j) v_s(\i)  + C_A v_s(\i) \Big)\, ds.
\end{equation}
Summing the previous equations over all $\i \in \V$,
interchanging the order of $\j$ and $\i$
and recalling that $\sum_{\i\in \V } p(\j,\i) = 1$, we get that
\begin{equation}\label{Ebound}
\bb{E}\crt{\norm{D_{t\wedge \sigma_A}}_1} \leq C_A\int_0^t \bb{E}\crt{\norm{D_{s\wedge \sigma_A}}_1} \, ds
\end{equation}
which means, by Gronwall inequality, that
\begin{equation}
\label{e:gb0}
\bb{E}\crt{\norm{D_{t\wedge \sigma_A}}_1} = 0.
\end{equation}
Since $\lim_{A \to \infty} \sigma_A = \infty$,
by Fatou's Lemma and \eqref{e:gb0} we obtain
\begin{equation}\label{fatou}
\bb{E}\crt{\norm{D_{t}}_1}
=    \bb{E}\crt{\liminf_{A\to \infty} \norm{D_{t\wedge \sigma_A}}_1}
\leq \liminf_{A\to \infty} \bb{E}\crt{\norm{D_{t\wedge \sigma_A}}_1}
=    0.
\qedhere
\end{equation}

\subsection{Solutions of martingale problem: sketch of proof of Proposition~\ref{prop:are_sols}} \label{s:are_sols}
\begin{proof}[Proof of Proposition~\ref{prop:are_sols}, sketch]
We write~$M^f_t:= f(\zeta_t)-f(\zeta_0) - \int_0^t \mathcal{L}^*f(\zeta_s)ds$.
From the fact that~$M^f_\cdot$ is a local martingale for every~$f \in \{f_x,f_{xy}:x,y \in \V\}$
and by repeating the computation in~\cite[Eq.(4.11), p.315]{KarShr98},
we obtain the cross-variation processes (see~\cite[Def.5.5, p.31]{KarShr98})  given by
\begin{equation}\label{eq:cross_var}
\langle M^{f_x}_\cdot,M^{f_y}_\cdot\rangle_t = \mathds{1}_{\{x=y\}}\cdot a\int_0^t  (\zeta_s(x))^\ell ds,\qquad t \ge 0.
\end{equation}

We now extend the probability space in which~$\zeta_\cdot$ is defined
to a space $(\hat{\Omega},\hat{\mathcal{F}},\hat{\P})$
where an auxiliary family of independent Brownian motions
$\{W^x_\cdot\}_{x \in \V}$ is defined.
To keep track of this extension,
we denote by~$\hat{\zeta}_\cdot$ and~$\hat{M}^f_\cdot$
the versions of~$\zeta_\cdot$ and~$M^f_\cdot$ in the larger space.
We define
\[
  B^x_t
  := \int_0^t \mathds{1}_{\{\hat{\zeta}_s(x)\neq 0\}}\cdot \frac{1}{\sqrt{a(\hat{\zeta}_s(x))^\ell}} \; d\hat{M}^{f_x}_s
    + \int_0^t \mathds{1}_{\{\hat{\zeta}_s(x) =0\}} \; dW^x_s.
\]
Using~\eqref{eq:cross_var}, it can then be proved that
$\{B^x_\cdot\}_{x \in \V}$ is a family of continuous local martingales
with cross-variation process given by
\[
  \langle B^x_\cdot,B^y_\cdot\rangle_t = t\cdot \mathds{1}_{\{x=y\}},\qquad t \ge 0.
\]
By Lévy's representation theorem~\cite[p.157]{KarShr98},
it follows that~$\{B^x_\cdot\}_{x \in \V}$
is a family of independent Brownian motions.
Since, almost surely,
\[
  \hat{M}^{f_x}_t
  = \int_0^t \left(a (\zeta_s(x))^\ell\right)^{1/2} dB^x_s, \quad t \ge 0
\]
and~$\mathcal{L}^*f_x(\zeta) = \Delta_p\zeta(x) - b(\zeta(x))^\kappa$,
we can, almost surely, rewrite the equation
\[
  \hat{\zeta}_t(x) = \hat{\zeta}_0(x) + \int_0^t \mathcal{L}^*f_x(\hat{\zeta}_s)ds + \hat{M}^{f_x}_t
\]
as
\[
  \hat{\zeta}_t(x)
  = \hat{\zeta}_0(x)
    + \int_0^t \left(\Delta_p\hat{\zeta}_s(x) - b(\hat{\zeta}_s(x))^\kappa\right) ds
    + \int_0^t \left(a(\hat{\zeta}_s(x))^\ell \right)^{1/2}dB^x_s.
  \qedhere
\]
\end{proof}

\subsection{Generator computations: proof of Lemma~\ref{lem:expression_generator}} \label{sec:app_gen}

\begin{proof}[Proof of Lemma~\ref{lem:expression_generator}]
Define
\[
  (\mathcal{S}^\m f)(\zeta) := R^\m(n\zeta) \cdot \left(f(\tfrac{1}{n}\Gamma^\m(n\zeta)) - f(\zeta) \right),
\]
for each~$\m \in \mathcal{M}$, so that
\begin{equation}\label{eq:main_xy}
  (\mathcal{L}^nf)(\zeta) = \sum_{\m\in \mathcal{M}}(\mathcal{S}^\m f)(\zeta).
\end{equation}

Note that
\[
  (\mathcal{S}^{(x,+)}f_x)(\zeta) = F^{n,+}(n\zeta(x))\cdot \frac{1}{n} = \frac{an}{2}\cdot (\zeta(x))^\ell
  - \min\left\{\frac{an}{2}\cdot (\zeta(x))^\ell;\;
    \frac{b}{2}\cdot (\zeta(x))^\kappa\right\}
\]
and similarly,
\[
  (\mathcal{S}^{(x,-)}f_x)(\zeta) =  \frac{an}{2}\cdot (\zeta(x))^\ell
  + \min\left\{\frac{an}{2}\cdot (\zeta(x))^\ell
    ;\; \frac{b}{2}\cdot (\zeta(x))^\kappa\right\}.
\]
Next, for~$x \neq y$,
\[
(\mathcal{S}^{(x,y)}f_x)(\zeta) = n\cdot \zeta(x)\cdot p(x,y)\cdot \left(-\frac{1}{n}\right) = -\zeta(x)\cdot p(x,y)
\]
and similarly,
\[
(\mathcal{S}^{(y,x)}f_x)(\zeta) = \zeta(y)\cdot p(y,x).
\]
Using these expressions in~\eqref{eq:main_xy} we obtain~\eqref{eq:gen_exp_Lnx}.

To prove~\eqref{eq:gen_exp_Qnx} we observe that
\begin{align*}
\mathcal{Q}^nf_x(\zeta) &= \sum_{\m \in \mathcal{M}}R^\m(n\zeta)\cdot \left(f_x(\tfrac{1}{n}\Gamma^\m(n\zeta))- f_x(\zeta)\right)^2\\
&=\left(R^{(x,+)}(n\zeta) + R^{(x,-)}(n\zeta) + \sum_{y \neq x}(R^{(x,y)}(n\zeta)+R^{(y,x)}(n\zeta))\right) \cdot \frac{1}{n^2}\\
&=a\cdot \zeta(x)^\ell + \frac{1}{n}\sum_{y \neq x}(p(y,x)\zeta(y) + p(x,y)\zeta(x)).
\end{align*}
We obtain~\eqref{eq:gen_exp_Lnxx} as a consequence of~\eqref{eq:gen_exp_Lnx}
and~\eqref{eq:gen_exp_Qnx} and the identity $\mathcal{L}^nf_{xx} = 2f_x\cdot \mathcal{L}^nf_x+ \mathcal{Q}^nf_{xx}$.
Finally, assume that~$x \neq y$, and let us prove~\eqref{eq:gen_exp_Lnxy}. We compute
\begin{align*}
  (\mathcal{S}^{(x,y)}f_{x,y})(\zeta)
  &= n\zeta(x)\cdot p(x,y)\cdot\left[\left(\zeta(x) - \frac{1}{n} \right)\cdot
    \left( \zeta(y) + \frac{1}{n}\right) -\zeta(x)\cdot \zeta(y) \right]\\
  &= p(x,y)\cdot \left(\zeta(x)^2 - \zeta(x)\cdot \zeta(y) - \frac{\zeta(x)}{n}\right).
\end{align*}
We then replace in~\eqref{eq:main_xy} the above, together with the analogous expression for~$\mathcal{S}^{(y,x)}$, to obtain
\begin{equation}
\label{eq:main_xy2}
\begin{split}
  (\mathcal{L}^nf_{x,y})(\zeta) = \sum_{\m\notin \{(x,y),(y,x)\}}(\mathcal{S}^\m f_{x,y})(\zeta)
  &+ p(x,y)\cdot \left(\zeta(x)^2 - \zeta(x)\cdot \zeta(y) - \frac{\zeta(x)}{n}\right) \\
  &\;\;+ p(y,x)\cdot \left(\zeta(y)^2 - \zeta(x)\cdot
    \zeta(y) - \frac{\zeta(y)}{n}\right).
\end{split}
\end{equation}

Next, for each~$z \in \V$, let~$\mathcal{M}_z := \{(z,+),(z,-)\}\cup\{(z,w): w \in \V\} \cup \{(w,z): w \in \V\}$. Note that
\[
  \m \in \mathcal{M}_x \backslash \mathcal{M}_y \quad \Longrightarrow \quad (\mathcal{S}^\m f_{x,y})(\zeta) = (\mathcal{S}^\m f_x)(\zeta)\cdot f_y(\zeta).
\]
This observation and a simple computation give
\begin{align*}
  (\mathcal{L}^n f_x)(\zeta)\cdot f_y(\zeta)
  &= \sum_{\m \in \mathcal{M}_x \backslash \mathcal{M}_y} (\mathcal{S}^\m f_{x,y})(\zeta) + (\mathcal{S}^{(x,y)}f_x)(\zeta)\cdot
    f_y(\zeta) +(\mathcal{S}^{(y,x)}f_x)(\zeta)\cdot  f_y(\zeta)\\
  &= \sum_{\m \in \mathcal{M}_x \backslash \mathcal{M}_y} (\mathcal{S}^\m f_{x,y})(\zeta) - p(x,y)\cdot
    \zeta(x)\cdot \zeta(y) + p(y,x)\cdot \zeta(y)^2
\end{align*}
and similarly,
\[
  f_x(\zeta)\cdot(\mathcal{L}^n f_y)(\zeta)
  = \sum_{\m \in \mathcal{M}_y \backslash \mathcal{M}_x} (\mathcal{S}^\m f_{x,y})(\zeta) - p(y,x)\cdot
  \zeta(x)\cdot \zeta(y) - p(x,y)\cdot \zeta(x)^2.
\]
Comparing these two expressions with~\eqref{eq:main_xy2} gives~\eqref{eq:gen_exp_Lnxy}.
\end{proof}

\bibliography{mybib_C.bib}
\bibliographystyle{plain}
\end{document}